# ON RECURSIVE ESTIMATION FOR TIME VARYING AUTOREGRESSIVE PROCESSES

BY ERIC MOULINES, PIERRE PRIOURET AND FRANÇOIS ROUEFF

*GET/Télécom Paris, CNRS LTCI, Université Paris VI
and GET/Télécom Paris, CNRS LTCI*

This paper focuses on recursive estimation of time varying autoregressive processes in a nonparametric setting. The stability of the model is revisited and uniform results are provided when the time-varying autoregressive parameters belong to appropriate smoothness classes. An adequate normalization for the correction term used in the recursive estimation procedure allows for very mild assumptions on the innovations distributions. The rate of convergence of the pointwise estimates is shown to be minimax in $\beta$-Lipschitz classes for $0 < \beta \leq 1$. For $1 < \beta \leq 2$, this property no longer holds. This can be seen by using an asymptotic expansion of the estimation error. A bias reduction method is then proposed for recovering the minimax rate.

**1. Introduction.** Suppose that we have real-valued observations $(X_{1,n}, X_{2,n}, \ldots, X_{n,n})$ from a time-varying autoregressive model (TVAR)

$$(1) \qquad X_{k,n} = \sum_{i=1}^{d} \theta_i((k-1)/n)X_{k-i,n} + \sigma(k/n)\varepsilon_{k,n}, \qquad k = 1, \ldots, n,$$

where

– $\{\varepsilon_{k,n}\}_{1 \leq k \leq n}$ is a triangular array of real-valued random variables referred to as the (normalized) *innovations*,
– $\boldsymbol{\theta}(t) := [\theta_1(t) \cdots \theta_d(t)]^T$, $t \in [0,1]$, is a $d$-dimensional vector referred to as the *local autoregression vector*,
– $\sigma(t)$, $t \in [0,1]$, is a nonnegative number referred to as the *local innovation standard deviation*.

This recurrence equation may be more compactly written as

$$(2) \qquad X_{k,n} = \boldsymbol{\theta}_{k-1,n}^T \mathbf{X}_{k-1,n} + \sigma_{k,n}\varepsilon_{k,n}, \qquad k = 1, \ldots, n,$$

where

$$\mathbf{X}_{k,n} := [X_{k,n} X_{k-1,n} \cdots X_{k-d+1,n}]^T,$$
$$\boldsymbol{\theta}_{k,n} := \boldsymbol{\theta}(k/n) = [\theta_1(k/n)\theta_2(k/n)\cdots\theta_d(k/n)]^T \quad \text{and} \quad \sigma_{k,n} = \sigma(k/n).$$

TVAR models have been used for modeling data whose spectral content varies along time (see, e.g., [12, 15, 22] for early references). TVAR models are also closely related to the general class of locally stationary processes (see [5, 6, 7] and Remark 3 below for definitions and properties).

In this paper we focus on the estimation of the functions $t \mapsto \boldsymbol{\theta}(t)$ [we leave aside $\sigma(t)$ for brevity] from the observations $\{\mathbf{X}_{0,n}, X_{k,n}, k \geq 1\}$ (here we add the initial conditions $\mathbf{X}_{0,n}$ in the observations set for convenience). This problem is reminiscent of nonparametric curve estimation on a fixed design, a problem which has received considerable attention in the literature. A natural approach consists of using a stationary method on short overlapping segments of the time series (see, e.g., [8]). An alternative approach, first investigated by [3] for first order TVAR models, consists of estimating the regression function recursively in time. More precisely, at a given time $t \in (0,1)$, only observations that have been observed *before* time $t$ are used in the definition of the estimator: $\hat{\boldsymbol{\theta}}_n(t) = \hat{\boldsymbol{\theta}}_n(t, \mathbf{X}_{0,n}, X_{1,n}, \ldots, X_{[nt],n})$, where $[x]$ denotes the integer part of $x$. This approach is useful when the observations must be processed *on line* (see, e.g., [18, 19, 21]). We focus in this contribution on the normalized least square algorithm (NLMS), which is a specific example of a recursive identification algorithm, defined as follows:

$$(3) \qquad \begin{aligned} &\hat{\boldsymbol{\theta}}_{0,n}(\mu) := 0, \\ &\hat{\boldsymbol{\theta}}_{k+1,n}(\mu) := \hat{\boldsymbol{\theta}}_{k,n}(\mu) + \mu(X_{k+1,n} - \hat{\boldsymbol{\theta}}_{k,n}^T(\mu)\mathbf{X}_{k,n})\frac{\mathbf{X}_{k,n}}{1 + \mu|\mathbf{X}_{k,n}|^2}, \end{aligned}$$

where $k$ goes from zero to $n-1$, $\mu$ is referred to as the step-size and $|\cdot|$ denotes the Euclidean norm. At each iteration of the algorithm, the parameter estimates are updated by moving in the direction of the gradient of the instantaneous estimate $(X_{k+1,n} - \boldsymbol{\theta}^T \mathbf{X}_{k,n})^2$ of the *local* mean square error $\mathbb{E}[(X_{k+1,n} - \boldsymbol{\theta}^T \mathbf{X}_{k,n})^2]$. The normalization $(1 + \mu|\mathbf{X}_{k,n}|^2)^{-1}$ is a safeguard again large values of the norm of the regression vector and allows for very mild assumptions on the innovations [see (A1) below] compared with the LMS, which typically requires much stronger assumptions (see [20]). Extensions to a more sophisticated iterative rule, for example, the so-called recursive least-square (RLS) algorithm, are currently under investigation.

We define a pointwise estimate of $t \mapsto \boldsymbol{\theta}(t)$ as a simple interpolation of $\hat{\boldsymbol{\theta}}_{k,n}(\mu), k = 1, \ldots, n$, that is,

$$(4) \qquad \hat{\boldsymbol{\theta}}_n(t; \mu) := \hat{\boldsymbol{\theta}}_{[tn],n}(\mu), \qquad t \in [0, 1], n \geq 1.$$

Observe that, for all $t \in [0, 1]$, $\hat{\boldsymbol{\theta}}_n(t; \mu)$ is a function of $\mathbf{X}_{0,n}$, $\{X_{l,n}, l = 1, \ldots, [tn]\}$ and $\mu$.

The paper is organized as follows. In Section 2 we introduce the assumptions and state the main results of this paper, namely, uniform risk bounds for $\hat{\boldsymbol{\theta}}_n$, a lower bound on the minimax $L^2$-risk and precise approximations of the risk for $\hat{\boldsymbol{\theta}}_n$. We also discuss a technique to correct the bias of the estimator. In Section 3 the basic results used in the paper for dealing with weak dependence are presented. In Section 4 a four-step proof of the uniform upper bound for the $L^p$-risk of $\hat{\boldsymbol{\theta}}_n$ is given. Section 5 then provides the proof of the minimax lower bound. In Section 6 further approximation results are given, from which we obtain the proofs of the risk approximations for $\hat{\boldsymbol{\theta}}_n$ stated in Section 2.

## 2. Main results.

The space of $m \times n$ matrices is embedded with the operator norm associated to the Euclidean norm, which we denote by

$$(5) \qquad |A| := \sup_{x \in \mathbb{R}^n, |x|=1} |Ax|.$$

Observe that for a row or column vector, its Euclidean norm coincides with its operator norm. For any random variable $\mathbf{Z}$ in a normed space $(\mathsf{Z}, |\cdot|)$, we denote by $\|\mathbf{Z}\|_p := (\mathbb{E}|\mathbf{Z}|^p)^{1/p}$ its $L^p$-norm. Throughout the paper, it is assumed that,

(A1) for all $n \geq 1$, the random variables $\{\varepsilon_{k,n}\}_{1 \leq k \leq n}$ are independent, have zero mean and unit variance and are independent of the initial conditions $\mathbf{X}_{0,n}$. In addition, $\sup_{n \geq 1} \|\mathbf{X}_{0,n}\|_q < \infty$ and $\varepsilon_q^\star := \sup_{1 \leq k \leq n} \|\varepsilon_{k,n}\|_q < \infty$, where the moment order $q \geq 2$ will be set depending on the context. The triangular array of random variables $\{X_{k,n}, 1 \leq k \leq n\}$ defined by (1) is parameterized by $(\boldsymbol{\theta}, \sigma)$. To keep track of the dependence on $(\boldsymbol{\theta}, \sigma)$, for all random variables $Z$ defined as a function of these variables, we shall adopt the notation convention $\mathbb{E}_{\boldsymbol{\theta}, \sigma}[Z]$ for the expectation of $Z$. In the case of a random element $\mathbf{Z}$ taking its values in the normed space $(\mathsf{Z}, |\cdot|)$, its $L^p$-norm will be denoted by

$$\|\mathbf{Z}\|_{p, \boldsymbol{\theta}, \sigma} := (\mathbb{E}_{\boldsymbol{\theta}, \sigma} |\mathbf{Z}|^p)^{1/p}.$$

A classical problem in nonparametric estimation is to derive uniform bounds for the pointwise $L^p$-risk $\|\hat{\boldsymbol{\theta}}_n(t; \mu) - \boldsymbol{\theta}(t)\|_{p, \boldsymbol{\theta}, \sigma}$ for $(\boldsymbol{\theta}, \sigma)$ in some appropriate

classes of functions. In the sequel we denote by

$$(6) \qquad \theta(z;t) := 1 - \sum_{j=1}^{d} \boldsymbol{\theta}_j(t) z^j, \qquad z \in \mathbb{C},$$

*the local time-varying autoregressive polynomial* associated with $\boldsymbol{\theta}$ at point $t$. The function classes that will be considered rely on two kinds of properties. First, the roots of the time-varying autoregressive polynomial associated with $\boldsymbol{\theta}$ are required to stay away from the unit disk. Second, as in function estimation from noisy data, $\boldsymbol{\theta}$ and $\sigma$ are supposed to be smooth in some appropriate sense. The first condition is unusual in the nonparametric function estimation setting and deserves some elaboration. Let us recall some elementary facts from the theory of autoregressive processes. The process $\{Z_k\}$ is an AR($d$) process with parameters $\boldsymbol{\vartheta} \in \mathbb{R}^d$ and $\varsigma^2 > 0$ if $\{Z_k\}$ is second-order stationary and satisfies the difference equation

$$(7) \qquad Z_k = \sum_{j=1}^{d} \boldsymbol{\vartheta}_j Z_{k-j} + \varepsilon_k, \qquad k \in \mathbb{Z},$$

where $\{\varepsilon_k\}_{k \in \mathbb{Z}}$ is a centered white noise with variance $\varsigma^2$. A sufficient and necessary condition for the existence of $\{Z_k\}$ is that the autoregressive polynomial $z \mapsto \vartheta(z) := 1 - \sum_{j=1}^{d} \boldsymbol{\vartheta}_j z^j$ does not vanish on the unit circle (see [4]). In this case, the stationary solution $\{Z_k\}$ to (7) is unique and there exists a sequence $\{\psi_k\}$, such that $\sum_l |\psi_l| < \infty$ and $Z_k = \sum_l \psi_l \varepsilon_{k-l}$ for all $k \in \mathbb{Z}$. Furthermore, the sequence $\{\psi_l\}$ is causal, that is, $\psi_l = 0$ for all $l < 0$, if and only if the function $z \mapsto \vartheta(z)$ does not vanish on the disk $|z| < \rho^{-1}$, for some $\rho < 1$. This motivates the following definitions in the time varying setting. For $\rho > 0$, we denote

$$(8) \quad \mathcal{S}(\rho) := \{\boldsymbol{\theta} : [0,1] \to \mathbb{R}^d, \theta(z;t) \neq 0 \text{ for all } |z| < \rho^{-1} \text{ and } t \in [0,1]\}.$$

Concerning the smoothness condition, different classes of functions can be considered. In the original paper by [6], it is assumed that the functions $t \mapsto \boldsymbol{\theta}(t)$ and $t \mapsto \sigma(t)$ are Lipschitzian. In this paper we consider a wider range of smoothness classes which are now introduced. For any $\beta \in (0,1]$, denote the $\beta$-Lipschitz semi-norm of a mapping $\mathbf{f} : [0,1] \mapsto \mathbb{R}^l$ by

$$|\mathbf{f}|_{\Lambda,\beta} = \sup_{t \neq s} \frac{|\mathbf{f}(t) - \mathbf{f}(s)|}{|t-s|^\beta}.$$

And define for $0 < L < \infty$ the $\beta$-Lipschitz ball

$$(9) \qquad \Lambda_l(\beta, L) := \left\{ \mathbf{f} : [0,1] \to \mathbb{R}^l, |\mathbf{f}|_{\Lambda,\beta} \leq L, \sup_{t \in [0,1]} |\mathbf{f}(t)| \leq L \right\}.$$

For all $\beta > 1$, the $\beta$-Lipschitz balls are classically generalized as follows. Let $k \in \mathbb{N}$ and $\alpha \in (0, 1]$ be uniquely defined by $\beta = k + \alpha$. Then we define

$$(10) \qquad \Lambda_l(\beta, L) := \left\{ \mathbf{f} : [0, 1] \to \mathbb{R}^l, |\mathbf{f}^{(k)}|_{\Lambda, \alpha} \leq L, \sup_{t \in [0,1]} |\mathbf{f}(t)| \leq L \right\},$$

where $\mathbf{f}^{(k)}$ is the derivative of order $k$.

For all $\beta > 0$, $L > 0$, $0 < \rho < 1$, and $0 < \sigma_- \leq \sigma_+ < \infty$, we define

$$\mathcal{C}(\beta, L, \rho, \sigma_-, \sigma_+) := \{(\boldsymbol{\theta}, \sigma) : \boldsymbol{\theta} \in \Lambda_d(\beta, L) \cap \mathcal{S}(\rho), \sigma : [0, 1] \to [\sigma_-, \sigma_+]\}.$$

We will simply write $\mathcal{C}$ whenever no confusion is possible. It is interesting to observe that, for particular choices of $L$ and $\rho$, $\mathcal{C}$ reduces to the more conventional smoothness class $\{(\boldsymbol{\theta}, \sigma) : \boldsymbol{\theta} \in \Lambda_d(\beta, L), \sigma : [0, 1] \to [\sigma_-, \sigma_+]\}$. This follows from the following lemma.

LEMMA 1. *For all positive $\rho$, we have*

$$(11) \qquad B(1/\sqrt{\rho^{-2} + \cdots + \rho^{-2d}}) \subseteq \mathcal{S}(\rho) \subseteq B((1 + \rho)^d - 1),$$

*where $B(a)$ is the sup-norm ball $\{\boldsymbol{\theta} : [0, 1] \to \mathbb{R}^d, \sup_{t \in [0,1]} |\boldsymbol{\theta}(t)| \leq a\}$.*

PROOF. Note that $\theta(t; 0) = 1$ for all $t \in [0, 1]$. Let $\lambda_1(t), \ldots, \lambda_d(t)$ be the reciprocals of the zeros of the polynomial $z \mapsto \theta(z; t)$. Hence, $\boldsymbol{\theta} \in \mathcal{S}(\rho)$ iff, for all $t \in [0, 1]$ and $k = 1, \ldots, d$, $|\lambda_k(t)| \leq \rho$. For all $t \in [0, 1]$ and $k = 1, \ldots, d$, using the Cauchy–Schwarz inequality, we have

$$1 = \left| \sum_{i=1}^d \theta_i(t) \lambda_k^{-i}(t) \right| \leq |\boldsymbol{\theta}(t)| \left( \sum_{i=1}^d |\lambda_k(t)|^{-2i} \right)^{1/2}.$$

Hence, the first inclusion in (11). We further have

$$\theta(z; t) = \prod_{k=1}^d (1 - \lambda_k(t) z).$$

The coefficients of $\boldsymbol{\theta}$ are then given by

$$(12) \qquad \theta_k(t) = (-1)^k \sum_{1 \leq i_1 < \cdots < i_k \leq d} \lambda_{i_1}(t) \cdots \lambda_{i_k}(t), \qquad k = 1, \ldots, d.$$

A simple computation then gives the second inclusion in (11). $\quad\square$

REMARK 1. In the case $d = 1$ considered in [3], $\mathcal{S}(\rho) = B(\rho)$.

We are now in a position where we can state the main results of this paper. We first provide a uniform upper bound on the pointwise $L^p$-risk of the NLMS estimator.

THEOREM 2. *Assume* (A1) *with* $q \geq 4$ *and let* $p \in [1, q/3]$. *Let* $\beta \in (0, 1]$, $L > 0$, $0 < \rho < 1$ *and* $0 < \sigma_- \leq \sigma_+$. *Then there exist* $M, \delta > 0$ *and* $\mu_0 > 0$ *such that, for all* $\mu \in (0, \mu_0]$, $n \geq 1$, $t \in (0, 1]$ *and* $(\boldsymbol{\theta}, \sigma) \in \mathcal{C}(\beta, L, \rho, \sigma_-, \sigma_+)$,

$$(13) \qquad \|\hat{\boldsymbol{\theta}}_n(t; \mu) - \boldsymbol{\theta}(t)\|_{p, \boldsymbol{\theta}, \sigma} \leq M(|\boldsymbol{\theta}(0)|(1 - \delta \mu)^{tn} + \sqrt{\mu} + (n \mu)^{-\beta}).$$

COROLLARY 3. *For all* $\eta \in (0, 1)$ *and* $\alpha > 0$, *there exists* $M > 0$ *such that, for all* $(\boldsymbol{\theta}, \sigma) \in \mathcal{C}(\beta, L, \rho, \sigma_-, \sigma_+)$,

$$\sup_{t \in [\eta, 1]} \|\hat{\boldsymbol{\theta}}_n[t; \alpha n^{-2\beta/(1+2\beta)}] - \boldsymbol{\theta}(t)\|_{p, \boldsymbol{\theta}, \sigma} \leq M n^{-\beta/(1+2\beta)}.$$

The upper bound in (13) has three terms. Anticipating what will be said in the proof section (Section 4), the first term $|\boldsymbol{\theta}(0)|(1 - \delta \mu)^{tn}$ reflects the forgetting of the initial error of the NLMS estimator. The second term is the so-called lag-noise term, which accounts for the fluctuation of the recursive estimator which would occur even if $t \mapsto \boldsymbol{\theta}(t)$ were constant. The third term controls the error involved by time evolution of $\boldsymbol{\theta}(t)$ and mainly relies on the smoothness exponent $\beta$. Corollary 3 is then obtained by choosing the step-size in order to minimize this upper bound. Observing that the first term is negligible for $t \geq \eta > 0$ and balancing the two remaining ones yields $\mu \propto n^{2\beta/(1+2\beta)}$. This corollary says that, for all $\beta \in (0, 1]$, under the $\beta$-Lipschitz assumption, for $t \in (0, 1]$ the $L^p$-risk of the NLMS estimator at point $t$ has rate $n^{-\beta/(1+2\beta)}$.

We now provide a lower bound on the $L^2$-risk for any estimator $\hat{\boldsymbol{\delta}}_n$ of $\boldsymbol{\theta}(t)$ computed from observations $\mathbf{X}_{0,n}, X_{1,n}, \ldots, X_{n,n}$. Let us stress that this lower bound is not restricted to recursive estimators, that is, we do not require $\hat{\boldsymbol{\delta}}_n$ to depend only on $\mathbf{X}_{0,n}, X_{1,n}, \ldots, X_{[nt],n}$. Denote by

$$\text{MSEM}_{\boldsymbol{\theta}, \sigma}(\hat{\boldsymbol{\delta}}_n, t) := \mathbb{E}_{\boldsymbol{\theta}, \sigma}[(\hat{\delta}_n - \boldsymbol{\theta}(t))(\hat{\delta}_n - \boldsymbol{\theta}(t))^T]$$

the mean square error matrix (MSEM) at $t \in [0, 1]$. Consider the following assumption:

(A2) For all $n \in \mathbb{N}$ and $1 \leq k \leq n$, $\varepsilon_{k,n}$ has an absolutely continuous density $p_{k,n}$ w.r.t. the Lebesgue measure whose derivative $\dot{p}_{k,n}$ satisfies

$$\mathbb{E}\left[\frac{\dot{p}_{k,n}}{p_{k,n}}(\varepsilon_{k,n})\right] = 0 \quad \text{and} \quad \mathcal{I}_\varepsilon := \sup_{1 \leq k \leq n} \mathbb{E}\left[\left(\frac{\dot{p}_{k,n}}{p_{k,n}}(\varepsilon_{k,n})\right)^2\right] < \infty.$$

We have the following result, the proof of which appears in Section 5.

THEOREM 4. *Assume* (A1) *with* $q = 2$ *and* (A2). *Let* $\beta > 0$, $L > 0$, $\rho \in (0, 1)$, $0 < \sigma_- \leq \sigma_+$. *Then there exists* $\alpha > 0$ *such that, for all* $n \geq 1$, $t \in [0, 1]$, *and for all estimators* $\hat{\boldsymbol{\delta}}_n := \hat{\delta}_n(\mathbf{X}_{0,n}, X_{1,n}, \ldots, X_{n,n}) \in \mathbb{R}^d$,

$$(14) \qquad \inf_{|\mathbf{u}|=1} \sup_{(\boldsymbol{\theta}, \sigma) \in \mathcal{C}} \mathbf{u}^T \text{MSEM}_{\boldsymbol{\theta}, \sigma}(\hat{\boldsymbol{\delta}}_n, t) \mathbf{u} \geq \alpha n^{-2\beta/(1+2\beta)},$$

*where* $\mathcal{C} := \mathcal{C}(\beta, L, \rho, \sigma_-, \sigma_+)$.

REMARK 2. Note that $\mathbf{u}^T \mathrm{MSEM}_{\boldsymbol{\theta},\sigma}(\hat{\boldsymbol{\delta}}_n, t)\mathbf{u}$ is the mean square error of $\mathbf{u}^T \hat{\boldsymbol{\delta}}_n$ for estimating $\mathbf{u}^T \boldsymbol{\theta}(t)$.

Corollary 3 and Theorem 4 show that, under (A1) with $q > 6$ and (A2), the $L^p$ error rate is minimax for $p \geq 2$ within the class $\mathcal{C}(\beta, L, \rho, \sigma_-, \sigma_+)$ if $\beta \in (0, 1]$. The question arises whether the upper bound derived in Theorem 2 generalizes for $\beta > 1$ in such a way that, as in Corollary 3, for an appropriate step-size $\mu(n)$, $\hat{\boldsymbol{\theta}}_n(t; \mu(n))$ achieves the rate of the lower bound derived in Theorem 4. It turns out that this is not the case, except in a very particular situation, namely, when $\boldsymbol{\theta}$ is the constant function. This may be shown by using precise approximations of the risk, completing the upper bound given in Theorem 2. Such approximations primarily rely on the fact that, as $n$ tend to infinity, the local sample $\mathbf{X}_{k,n}$ of the TVAR process approximately has the same second-order statistics as the stationary AR($d$) process with parameter $(\boldsymbol{\theta}(k/n), \sigma(k/n))$. In the following we provide a precise statement of this fact and then state the approximations of the risk. For this purpose, we need to introduce further notation. For $\beta > 0$, $L > 0$, $\rho \in (0, 1)$ and $0 < \sigma_- \leq \sigma_+$, we let

$$(15) \qquad \mathcal{C}^{\star}(\beta, L, \rho, \sigma_-, \sigma_+) := \{(\boldsymbol{\theta}, \sigma) \in \mathcal{C}(\beta, L, \rho, \sigma_-, \sigma_+) : \sigma \in \Lambda_1(\beta, \sigma_+)\}.$$

We use the shorthand notation $\mathcal{C}^{\star}$ when no confusion is possible. The obvious relation

$$\mathcal{C}(\beta, L, \rho, \sigma_-, \sigma_+) \supseteq \mathcal{C}^{\star}(\beta, L, \rho, \sigma_-, \sigma_+) \supseteq \mathcal{C}(\beta, L, \rho, \sigma_+, \sigma_+)$$

implies that Theorem 2 and Theorem 4 are still valid when replacing $\mathcal{C}$ by $\mathcal{C}^{\star}$. Following the formula for spectral densities of stationary AR($d$) processes, we respectively denote

$$(16) \qquad f(\lambda; t, \boldsymbol{\theta}, \sigma) := \frac{\sigma^2(t)}{2\pi} |\theta(e^{i\lambda}; t)|^{-2}, \qquad -\pi \leq \lambda \leq \pi,$$

$$(17) \qquad [\Sigma(t, \boldsymbol{\theta}, \sigma)]_{k,l} := \int_{-\pi}^{\pi} e^{i\lambda(k-l)} f(\lambda; t, \boldsymbol{\theta}, \sigma)\, d\lambda, \qquad 1 \leq k, l \leq d,$$

*the local spectral density function* and *the local $d$-dimensional covariance matrix* associated with $(\boldsymbol{\theta}, \sigma)$ at point $t$. The covariance matrix $\mathbb{E}_{\boldsymbol{\theta},\sigma}[\mathbf{X}_{k,n}\mathbf{X}_{k,n}^T]$ can be approximated by the local covariance matrix at point $k/n$ as follows.

PROPOSITION 5. *Assume* (A1) *with* $q \geq 2$. *Let* $\beta \in (0, 1]$, $L > 0$, $0 < \rho < \tau < 1$, *and* $0 < \sigma_- \leq \sigma_+$. *Then there exists* $M > 0$ *such that, for all* $1 \leq k \leq n$ *and* $(\boldsymbol{\theta}, \sigma) \in \mathcal{C}^{\star}(\beta, L, \rho, \sigma_-, \sigma_+)$,

$$|\mathbb{E}_{\boldsymbol{\theta},\sigma}[\mathbf{X}_{k,n}\mathbf{X}_{k,n}^T] - \Sigma(k/n, \boldsymbol{\theta}, \sigma)| \leq M(\tau^k |\mathbb{E}[\mathbf{X}_{0,n}\mathbf{X}_{0,n}^T] - \Sigma(0, \boldsymbol{\theta}, \sigma)| + n^{-\beta}).$$

REMARK 3. This approximation result can serve as an illustration of how the TVAR process fits into the locally stationary setting introduced by Dahlhaus [5]. Observe that if $\beta = 1$ and $\mathbb{E}[\mathbf{X}_{0,n}\mathbf{X}_{0,n}^T] = \Sigma(0,\boldsymbol{\theta},\sigma)$, the rate for the approximation error between the local sample covariance matrix $\mathbb{E}_{\boldsymbol{\theta},\sigma}[\mathbf{X}_{k,n}\mathbf{X}_{k,n}^T]$ and its local *stationary* approximation $\Sigma(k/n,\boldsymbol{\theta},\sigma)$ is $n^{-1}$, which coincides with the approximation rate required in the locally stationary setting introduced in [5]. Precise conditions upon which a TVAR process is locally stationary are given in [6].

We obtain the following computational approximation of the pointwise MSEM for $\hat{\boldsymbol{\theta}}_n$.

THEOREM 6. *Assume* (A1) *with $q > 11$. Let $\beta \in (0,1]$, $L > 0$, $\rho < 1$, and $0 < \sigma_- \leq \sigma_+$ and let $(\boldsymbol{\theta},\sigma) \in \mathcal{C}^\star(\beta,L,\rho,\sigma_-,\sigma_+)$. Let $t \in (0,1]$ and assume that there exists $\boldsymbol{\theta}_{t,\beta} \in \mathbb{R}^d$, $L' > 0$ and $\beta' > \beta$ such that, for all $u \in [0,t]$,*

$$(18) \qquad |\boldsymbol{\theta}(u) - \boldsymbol{\theta}(t) - \boldsymbol{\theta}_{t,\beta}(t-u)^\beta| \leq L'(t-u)^{\beta'}.$$

*Then there exist $M > 0$ and $\mu_0 > 0$ such that, for all $\mu \in (0,\mu_0]$ and $n \geq 1$,*

$$
(19) \quad
\begin{aligned}
&\left| \mathrm{MSEM}_{\boldsymbol{\theta},\sigma}\left( \hat{\boldsymbol{\theta}}_n(t;\mu) - \frac{\Gamma(\beta+1)}{(\mu n)^\beta}\Sigma^{-\beta}(t,\boldsymbol{\theta},\sigma)\boldsymbol{\theta}_{t,\beta} \right) - \mu \frac{\sigma^2(t)}{2}I \right|^{1/2} \\
&\qquad \leq M(\sqrt{\mu}(\sqrt{\mu} + (\mu n)^{-\beta/2}) + (\mu n)^{-\beta}((\mu n)^{-\beta} + (\mu n)^{\beta-\beta'} + \sqrt{\mu})),
\end{aligned}
$$

*where $\Gamma$ is the usual Gamma function and $I$ is the $d \times d$ identity matrix.*

REMARK 4. Let $\alpha \in \mathbb{R}$. The $\alpha$-fractional power of a diagonal matrix $D$ with positive diagonal entries is the diagonal matrix $D^\alpha$ obtained by raising the diagonal entries to the power $\alpha$. The $\alpha$-fractional power of a symmetric positive-definite matrix $A = UDU^T$, where $U$ is unitary and $D$ is diagonal, is then defined by $A^\alpha = UD^\alpha U^T$.

Using (19), as $(\mu + (\mu n)^{-1}) \to 0$, we have the following asymptotic approximation of the MSEM:

$$\mathrm{MSEM}_{\boldsymbol{\theta},\sigma}(\hat{\boldsymbol{\theta}}_n(t;\mu)) - \Gamma(\beta+1)(\mu n)^{-\beta}\Sigma^{-\beta}(t,\boldsymbol{\theta},\sigma)\boldsymbol{\theta}_{t,\beta} = \mu\frac{\sigma^2(t)}{2}I(1+o(1)).$$

If $\beta \leq 1$, this gives the leading term of an asymptotic expansion of the MSEM, which allows one to compare the performance of $\hat{\boldsymbol{\theta}}_n$ with other estimators achieving the minimax rate. The deterministic correction $\Gamma(\beta+1)(\mu n)^{-\beta}\Sigma^{-\beta}(t,\boldsymbol{\theta},\sigma)\boldsymbol{\theta}_{t,\beta}$ can be interpreted as the main term of the bias and the term $\mu\frac{\sigma^2(t)}{2}I$ as the main term of the covariance matrix.

REMARK 5. An essential ingredient for proving Theorem 6 is to approximate the expectation of the term $\mathbf{X}_{k,n}^T \mathbf{X}_{k,n}/(1 + \mu|\mathbf{X}_{k,n}|^2)$ appearing in (3). Roughly speaking, for small enough $\mu$, a good approximation is $\mathbb{E}_{\boldsymbol{\theta},\sigma}[\mathbf{X}_{k,n}^T \mathbf{X}_{k,n}]$, which itself is well approximated by $\Sigma(k/n, \boldsymbol{\theta}, \sigma)$ by using Proposition 5. If one replaces the normalization factor $1 + \mu|\mathbf{X}_{k,n}|^2$ by the more classical (in the stochastic tracking literature) $1 + |\mathbf{X}_{k,n}|^2$, the computation of the deterministic approximation would be much more involved as the normalization does not reduce to one as $\mu$ tend to zero.

REMARK 6. Equation (18) is, for instance, valid if $\boldsymbol{\theta}$ behaves as a sum of nonentire power laws at the left of point $t$, say, $\boldsymbol{\theta}(u) = \boldsymbol{\theta}_{t,1/2}\sqrt{t-u} + O(t-u)$. Another simple case consists in assuming that $(\boldsymbol{\theta}, \sigma) \in \mathcal{C}^\star(\beta', L, \rho, \sigma_-, \sigma_+)$ for $\beta' > 1$. Then (18) is obtained with $\beta = 1$ and $\boldsymbol{\theta}_{t,\beta} = -\dot{\boldsymbol{\theta}}(t)$ by using a first-order Taylor expansion. Hence, in this case, the main terms of the MSEM are of order $\mu + (\mu n)^{-2}$ unless $\dot{\boldsymbol{\theta}}(t) = 0$, in which case the deterministic correction in the MSEM vanishes. This implies that the estimator $\hat{\boldsymbol{\theta}}_n(t; \mu)$ cannot achieve the minimax rate obtained in Theorem 4 in the class $(\boldsymbol{\theta}, \sigma) \in \mathcal{C}^\star(\beta', L, \rho, \sigma_-, \sigma_+)$ for $\beta' > 1$ unless $\boldsymbol{\theta}$ is a constant function.

If the smoothness exponent belongs to $(1, 2]$, the following result applies.

THEOREM 7. Assume (A1) with $q > 4$ and let $p \in [1, q/4)$. Let $\beta \in (1, 2]$, $L > 0$, $\rho \in (0, 1)$, and $0 < \sigma_- \leq \sigma_+$. Then, for all $\eta \in (0, 1)$, there exist $M > 0$ and $\mu_0 > 0$ such that, for all $(\boldsymbol{\theta}, \sigma) \in \mathcal{C}^\star(\beta, L, \rho, \sigma_-, \sigma_+)$, $t \in [\eta, 1]$, $n \geq 1$ and $\mu \in (0, \mu_0]$,

$$\|\hat{\boldsymbol{\theta}}_n(t; \mu) - \boldsymbol{\theta}(t) + (\mu n)^{-1}\Sigma^{-1}(t, \boldsymbol{\theta}, \sigma)\dot{\boldsymbol{\theta}}(t)\|_{p,\boldsymbol{\theta},\sigma} \leq M(\sqrt{\mu} + (\mu n)^{-\beta} + (\mu n)^{-2}).$$

Applying a technique inspired by the so-called Romberg method in numerical analysis (see, e.g., [2]), we are now able to propose a recursive estimator which achieves the minimax rates for $\beta \in (1, 2]$. This estimator is obtained by combining the recursive estimators $\hat{\boldsymbol{\theta}}(t; \cdot)$ associated to two different stepsizes. More precisely, let

$$\tilde{\boldsymbol{\theta}}_n(t; \mu, \gamma) := \frac{1}{1 - \gamma}(\hat{\boldsymbol{\theta}}_n(t; \mu) - \gamma\hat{\boldsymbol{\theta}}_n(t; \gamma\mu)),$$

where $\gamma \in (0, 1)$. We obtain the following result.

THEOREM 8. Assume (A1) with $q > 4$ and let $p \in [1, q/4)$. Let $\beta \in (1, 2]$, $L > 0$, $\rho < 1$, and $0 < \sigma_- \leq \sigma_+$. For all $\eta \in (0, 1)$, there exist $M > 0$ and $\mu_0 > 0$ such that, for all $\gamma \in (0, 1)$, $(\boldsymbol{\theta}, \sigma) \in \mathcal{C}^\star(\beta, L, \rho, \sigma_-, \sigma_+)$, $n \geq 1$ and

$\mu \in (0, \mu_0]$,

$$\sup_{t \in [\eta, 1]} \|\tilde{\boldsymbol{\theta}}_n(t; \mu, \gamma) - \boldsymbol{\theta}(t)\|_{p, \boldsymbol{\theta}, \sigma}$$

(20)

$$\leq M \frac{1 + \gamma}{\gamma(1 - \gamma)} (\sqrt{\mu} + (\mu n)^{-\beta} + (\mu n)^{-2}).$$

PROOF. Let $\eta \in (0, 1)$, $\gamma \in (0, 1)$ and $t \in [\eta, 1]$. One easily checks that

$$\tilde{\boldsymbol{\theta}}_n(t; \mu, \gamma) - \boldsymbol{\theta}(t)$$

$$= (1 - \gamma)^{-1}(\hat{\boldsymbol{\theta}}_n(t; \mu) - \boldsymbol{\theta}(t) + (\mu n)^{-1} \Sigma^{-1}(t, \boldsymbol{\theta}, \sigma) \dot{\boldsymbol{\theta}}(t)$$

$$- \gamma(\hat{\boldsymbol{\theta}}_n(t; \gamma\mu) - \boldsymbol{\theta}(t) + (\gamma\mu n)^{-1} \Sigma^{-1}(t, \boldsymbol{\theta}, \sigma) \dot{\boldsymbol{\theta}}(t))).$$

The Minkowski inequality, Theorem 7 and the bounds $\sqrt{\gamma} < \gamma^{-\beta} \leq \gamma^{-2}$ yield (20). □

COROLLARY 9. *For all* $\gamma, \eta \in (0, 1)$, $\beta \in (0, 2]$ *and* $\alpha > 0$, *there exists* $M > 0$ *such that, for all* $(\boldsymbol{\theta}, \sigma) \in \mathcal{C}^\star(\beta, L, \rho, \sigma_-, \sigma_+)$,

$$\sup_{t \in [\eta, 1]} \|\tilde{\boldsymbol{\theta}}_n(t; \alpha n^{-2\beta/(1+2\beta)}, \gamma) - \boldsymbol{\theta}(t)\|_{p, \boldsymbol{\theta}, \sigma} \leq M \frac{1 + \gamma}{\gamma(1 - \gamma)} n^{-\beta/(1+2\beta)}.$$

**3. Exponential stability of inhomogeneous difference equations.** Let us consider a sequence $\{\mathbf{Z}_k, k \geq 0\}$ of random vectors satisfying the inhomogeneous difference equation

(21)                    $\mathbf{Z}_k = A_k \mathbf{Z}_{k-1} + B_k \mathbf{U}_k$,       $k \geq 1$,

where $\{\mathbf{U}_k, k \geq 1\}$ is a sequence of independent random vectors and $A = \{A_k, k \geq 1\}, B = \{B_k, k \geq 1\}$ are two sequences of deterministic matrices with appropriate dimensions. The pair $(A, B)$ is said to be *exponentially stable* if there exist constants $C > 0$ and $\rho \in (0, 1)$ such that

(22)  $\sup_{k \geq 0} \left| \prod_{l=k+1}^{k+m} A_l \right| \leq C \rho^m$      for all $m > 0$   and   $B^\star := \sup_{k \geq 1} |B_k| < \infty$,

with the convention that $\prod_{l=k}^{k+m} A_l := A_{k+m} A_{k+m-1} \cdots A_k$. We clearly have the following:

PROPOSITION 10. *Let* $p \in [1, \infty]$. *Suppose that* $(A, B)$ *is exponentially stable. Then there exists a positive constant* $M$ *depending only on* $C, \rho$ *and* $B^\star$ *such that*

$$\|\mathbf{Z}_k\|_p \leq M(\mathbf{U}^\star p + \|\mathbf{Z}_0\|_p),       k \in \mathbb{N},$$

*where* $\mathbf{U}^\star p := \sup_{k \geq 1} \|\mathbf{U}_k\|_p$.

Exponential stability implies exponential forgetting of the initial condition in the following setting. Let $(\mathsf{E}, |\cdot|_\mathsf{E})$ and $(\mathsf{F}, |\cdot|_\mathsf{F})$ be two normed spaces, let $m$ be a positive integer and $p$ a nonnegative real number. We denote by $\mathrm{Li}(\mathsf{E}, m, \mathsf{F}; p)$ the linear space of mappings $\phi: \mathsf{E}^m \to \mathsf{F}$ for which there exists $\lambda_1 \geq 0$ such that, for all $(x_1, \ldots, x_m)$ and $(y_1, \ldots, y_m) \in \mathsf{E}^m$,

$$
\begin{aligned}
(23) \quad & |\phi(x_1, \ldots, x_m) - \phi(y_1, \ldots, y_m)|_\mathsf{F} \\
& \leq \lambda_1(|x_1 - y_1|_\mathsf{E} + \cdots + |x_m - y_m|_\mathsf{E}) \\
& \quad \times (1 + |x_1|_\mathsf{E}^p + |y_1|_\mathsf{E}^p + \cdots + |x_m|_\mathsf{E}^p + |y_m|_\mathsf{E}^p).
\end{aligned}
$$

This implies that there exists $\lambda_2 \geq 0$ such that, for all $(x_1, \ldots, x_m) \in \mathsf{E}^m$,

$$
(24) \quad |\phi(x_1, \ldots, x_m)|_\mathsf{F} \leq \lambda_2(1 + |x_1|_\mathsf{E}^{p+1} + \cdots + |x_m|_\mathsf{E}^{p+1}).
$$

Denote by $|\phi|_{\mathrm{Li}(p)}$ the smallest $\lambda$ satisfying (23) and (24). $\mathrm{Li}(\mathsf{E}, m, \mathsf{F}; p)$ is called the linear space of $p$-weighted Lipschitz mappings and $|\phi|_{\mathrm{Li}(p)}$ the $p$-weighted Lipschitz norm of $\phi$. We now state the exponential forgetting property. In the sequel, for any integrable r.v. $Z$, $\mathbb{E}^\mathcal{F}[Z]$ denotes the conditional expectation of $Z$ given the $\sigma$-field $\mathcal{F}$ and inequalities involving random variables will be meant in the *almost sure* sense.

PROPOSITION 11. *Let $p \geq 0$. Assume that $\mathbf{U}_{p+1}^\star$ is finite and that (22) is satisfied for some $C > 0$ and $\rho < 1$. Let $\phi \in \mathrm{Li}(\mathbb{R}^d, m, \mathbb{R}; p)$, where $m$ is a positive integer. For all $k \geq 0$, let us denote $\mathcal{F}_k := \sigma(\mathbf{Z}_0, \mathbf{U}_l, 1 \leq l \leq k)$. Then there exist constants $C_1$ and $C_2$ (depending only on $C, \rho, \mathbf{U}_{p+1}^\star, B^\star$ and $p$) such that, for all $0 \leq l \leq k$ and $0 = j_1 < \cdots < j_m$,*

$$
\begin{aligned}
& |\mathbb{E}^{\mathcal{F}_k}[\mathbf{\Phi}_{k+r}]| \leq C_1 m |\phi|_{\mathrm{Li}(p)}(1 + \rho^{r(p+1)}|\mathbf{Z}_k|^{p+1}), \\
(25) \quad & |\mathbb{E}^{\mathcal{F}_k}[\mathbf{\Phi}_{k+r}] - \mathbb{E}^{\mathcal{F}_l}[\mathbf{\Phi}_{k+r}]| \leq C_2 m \rho^r |\phi|_{\mathrm{Li}(p)}(1 + |\mathbf{Z}_k|^{p+1} + \mathbb{E}^{\mathcal{F}_l}|\mathbf{Z}_k|^{p+1}), \\
& |\mathbb{E}^{\mathcal{F}_k}[\mathbf{\Phi}_{k+r}] - \mathbb{E}[\mathbf{\Phi}_{k+r}]| \leq C_2 m \rho^r |\phi|_{\mathrm{Li}(p)}(1 + |\mathbf{Z}_k|^{p+1} + \mathbb{E}[|\mathbf{Z}_k|^{p+1}]),
\end{aligned}
$$

*where $\mathbf{\Phi}_j = \phi(\mathbf{Z}_{j+j_1}, \ldots, \mathbf{Z}_{j+j_m})$ for all $j \geq 0$.*

PROOF. Define $H_{k,r}$ as the $\mathbb{R}^d \to \mathbb{R}$ function mapping $x \in \mathbb{R}^d$ to

$$
\mathbb{E}[\phi(h_{k,r+j_1}(\mathbf{x}, \mathbf{U}_{k+1}, \ldots, \mathbf{U}_{k+r+j_1}), \ldots, h_{k,r+j_m}(\mathbf{x}, \mathbf{U}_{k+1}, \ldots, \mathbf{U}_{k+r+j_m}))],
$$

where, for all $i, j \in \mathbb{N}$ and $(\mathbf{x}, \mathbf{u}_1, \ldots, \mathbf{u}_j) \in \mathbb{R}^d \times \mathbb{R}^{qj}$,

$$
h_{i,j}(\mathbf{x}, \mathbf{u}_1, \ldots, \mathbf{u}_j) = \alpha(i+j, i)\mathbf{x} + \sum_{k=1}^{j} \alpha(i+j, i+k) B_{i+k} \mathbf{u}_k
$$

with $\alpha(l, l) = I$ and $\alpha(l+m, l) := \prod_{k=l+1}^{l+m} A_k, l \geq 0, m \geq 1$.

We thus have, for all $k, r \geq 0$,

$$\mathbb{E}^{\mathcal{F}_k}[\mathbf{\Phi}_{k+r}] = H_{k,r}(\mathbf{Z}_k). \tag{26}$$

Using (22), observe that, for all $i, j \in \mathbb{N}$ and $(\mathbf{x}, \mathbf{y}, \mathbf{u}_1, \ldots, \mathbf{u}_j) \in \mathbb{R}^{2d} \times \mathbb{R}^{qj}$,

$$|h_{i,j}(\mathbf{x}, \mathbf{u}_1, \ldots, \mathbf{u}_j) - h_{i,j}(\mathbf{y}, \mathbf{u}_1, \ldots, \mathbf{u}_j)| \leq C\rho^j |\mathbf{x} - \mathbf{y}|,$$

$$|h_{i,j}(\mathbf{x}, \mathbf{u}_1, \ldots, \mathbf{u}_j)| \leq C\left(\rho^j |\mathbf{x}| + \sum_{k=1}^{j} \rho^{j-k} B^\star |\mathbf{u}_k|\right).$$

Using these bounds with the Minkowski inequality, the definition of $H_{k,r}$ above and the assumptions on $\phi$, we easily obtain, for all $(\mathbf{x}, \mathbf{y}) \in \mathbb{R}^{2d}$,

$$|H_{k,r}(\mathbf{x}) - H_{k,r}(\mathbf{y})| \leq c_1 |\phi|_{\mathrm{Li}(p)}\left(\sum_{i=1}^{m} \rho^{r+j_i}\right)|\mathbf{x} - \mathbf{y}|$$

$$\times \left[1 + \sum_{i=1}^{m} \rho^{(r+j_i)p}(|\mathbf{x}|^p + |\mathbf{y}|^p)\right.$$

$$\left. + 2\sum_{i=1}^{m}\left(\sum_{k=1}^{r+j_i} \rho^{r+j_i-k} B^\star \mathbf{U}_p^\star\right)^p\right],$$

$$|H_{k,r}(\mathbf{x})| \leq c_1 |\phi|_{\mathrm{Li}(p)}$$

$$\times \left[1 + \sum_{i=1}^{m} \rho^{(r+j_i)(p+1)}|\mathbf{x}|^{p+1}\right.$$

$$\left. + \sum_{i=1}^{m}\left(\sum_{k=1}^{r+j_i} \rho^{r+j_i-k} B^\star \mathbf{U}_{p+1}^\star\right)^{p+1}\right],$$

where $c_1$ is a constant depending only on $C$, $\rho$ and $p$. Observing that $\sum_{k=1}^{j} \rho^{j-k} \leq 1/(1-\rho)$ and $\sum_{i=1}^{m} \rho^{j_i \alpha} \leq 1/(1-\rho^\alpha)$ for all $j \geq 1$, $0 = j_1 < \cdots < j_m$ and $\alpha > 0$, we obtain

$$|H_{k,r}(\mathbf{x})| \leq C_1 |\phi|_{\mathrm{Li}(p)} m(1 + \rho^{r(p+1)}|\mathbf{x}|^{p+1}), \tag{27}$$

$$|H_{k,r}(\mathbf{x}) - H_{k,r}(\mathbf{y})| \leq C_1 m \rho^r |\phi|_{\mathrm{Li}(p)}|x - y|(1 + \rho^{rp}(|\mathbf{x}|^p + |\mathbf{y}|^p)), \tag{28}$$

where $C_1$ is a constant depending only on $B^\star, \mathbf{U}_{p+1}^\star$, $C$, $\rho$ and $p$. Equation (25) follows from (26) and (27). Observe now that, for any probability measure $\zeta$ of $\mathbb{R}^d$,

$$\left|H_{k,r}(\mathbf{x}) - \int H_{k,r}(\mathbf{y})\zeta(d\mathbf{y})\right| \leq \int |H_{k,r}(\mathbf{x}) - H_{k,r}(\mathbf{y})|\zeta(d\mathbf{y})$$

$$\leq C_1 m \rho^r |\phi|_{\mathrm{Li}(p)}\int |\mathbf{x} - \mathbf{y}|(1 + |\mathbf{x}|^p + |\mathbf{y}|^p)\zeta(d\mathbf{y}),$$

where we have used (28). The two last bounds of the proposition follow by using (26), $|\mathbf{x} - \mathbf{y}| \leq |\mathbf{x}| + |\mathbf{y}|$ and by choosing $\zeta$ respectively equal to the conditional distribution of $\mathbf{Z}_k$ given $(\mathbf{Z}_0, \mathbf{U}_1, \ldots, \mathbf{U}_l)$ and to the distribution of $\mathbf{Z}_k$. $\quad\square$

**4. Proof of Theorem 2.** We denote by $\mathbb{M}_d$ and $\mathbb{M}_d^+$ the space of $d \times d$ real matrices and the subspace of positive semi-definite symmetric matrices, respectively. For all $A \in \mathbb{M}_d$, we let $\lambda_{\min}(A)$, $\lambda_{\max}(A)$ and $|\lambda|_{\max}(A)$ denote the minimum eigenvalue, the maximum eigenvalue and the spectral radius of the matrix $A$, respectively, that is, $\lambda_{\min}(A) := \min(\mathrm{sp}(A))$, $\lambda_{\max}(A) := \max(\mathrm{sp}(A))$ and $|\lambda|_{\max}(A) := \max|\mathrm{sp}(A)|$, where $\mathrm{sp}(A)$ denotes the set of eigenvalues of $A$. The proof is derived in five steps and relies on intermediary results which will be repeatedly used throughout the paper.

*Step* 1. *Exponential Stability of the TVAR model.* We have the following.

LEMMA 12. *Let $0 < \rho_0 < \rho$ and $L$ be a positive constant. Then there exists $C > 0$ such that, for all $A \in \mathbb{M}_d$ with $|A| \leq L$ and $|\lambda|_{\max}(A) \leq \rho_0$, and for all $k \in \mathbb{N}$, $|A^k| \leq C\rho^k$.*

PROOF. We apply Theorem VII.1.10 of [10]. Let $\boldsymbol{\gamma} = \{z \in \mathbb{C} : |z| = \rho\}$. Then for any $A$ such that $|\lambda|_{\max}(A) \leq \rho_0$ and for all $k \in \mathbb{N}$,

$$A^k = \frac{1}{2\pi i} \int_{\boldsymbol{\gamma}} z^k (zI - A)^{-1} \, dz.$$

Hence, putting $z = \rho e^{iy}$, $y \in (-\pi, \pi)$, we have

$$|A^k| \leq \frac{\rho^k}{2\pi} \int_{-\pi}^{\pi} |(I - A/(\rho e^{iy}))^{-1}| \, dy.$$

Let $\mathsf{G} := \{A \in \mathbb{M}_d : |A| \leq L, |\lambda|_{\max}(A) \leq \rho_0\}$. $\mathsf{G}$ is a compact set and $(z, A) \mapsto |(I - A/z)^{-1}|$ is continuous over $\boldsymbol{\gamma} \times \mathsf{G}$ so that it is uniformly bounded. The proof follows. $\quad\square$

For all $t \in [0, 1]$ and for all $\boldsymbol{\theta} = (\theta_1, \ldots, \theta_d) : [0, 1] \to \mathbb{R}^d$, let $\Theta(t, \boldsymbol{\theta})$ denote the companion matrix defined by

$$(29) \qquad \Theta(t, \boldsymbol{\theta}) = \begin{bmatrix} \theta_1(t) & \theta_2(t) & \ldots & \ldots & \theta_d(t) \\ 1 & 0 & \ldots & \ldots & 0 \\ 0 & 1 & 0 & \ldots & 0 \\ \vdots & & & & \vdots \\ 0 & \ldots & 0 & 1 & 0 \end{bmatrix},$$

whose eigenvalues are the reciprocals of the zeros of the autoregressive polynomial. Using this notation, (2) can be rewritten as

$$(30) \qquad \mathbf{X}_{k+1,n} = \Theta(k/n, \boldsymbol{\theta})\mathbf{X}_{k,n} + \boldsymbol{\sigma}_{k+1,n}\varepsilon_{k+1,n}, \quad 0 \le k \le n-1,$$

where $\boldsymbol{\sigma}_{k,n} = [\sigma(k/n)\, 0\cdots 0]^T$.

PROPOSITION 13. *Let $\beta \in (0,1]$, $L > 0$ and $0 < \rho < \tau < 1$. Then there exists a constant $M > 0$ such that, for all $\boldsymbol{\theta} \in \Lambda_d(\beta, L) \cap \mathcal{S}(\rho)$ and $0 \le k < k + m \le n$,*

$$(31) \qquad \left| \prod_{l=k+1}^{k+m} \Theta(l/n, \boldsymbol{\theta}) \right| \le M\tau^m.$$

PROOF. For notational convenience, we use $\Lambda$ and $\mathcal{S}$ as shorthand notation for $\Lambda_d(\beta, L)$ and $\mathcal{S}(\rho)$. First note that $\Theta^\star = \sup_{\boldsymbol{\theta} \in \Lambda \cap \mathcal{S}} \sup_{t \in [0,1]} |\Theta(t, \boldsymbol{\theta})|$ is finite. For any square matrices $A_1, \dots, A_r$ we have

$$(32) \qquad \begin{aligned} \prod_{k=1}^{r} A_k &= A_1^r + (A_r - A_1)A_1^{r-1} \\ &\quad + A_r(A_{r-1} - A_1)A_1^{r-2} + \cdots + \prod_{k=3}^{r} A_k(A_2 - A_1)A_1. \end{aligned}$$

By applying this decomposition to

$$(33) \qquad \begin{aligned} \beta_n(k, i; \boldsymbol{\theta}) &:= \prod_{j=i+1}^{k} \Theta(j/n, \boldsymbol{\theta}), \qquad 0 \le i < k \le n, \\ \beta_n(i, i; \boldsymbol{\theta}) &= I, \qquad\qquad\qquad\quad 0 \le i \le n, \end{aligned}$$

we have, for all $0 \le k < k + q \le n$,

$$\begin{aligned} |\beta_n(k+q, k; \boldsymbol{\theta})| &\le |\Theta((k+1)/n, \boldsymbol{\theta})^q| \\ &\quad + q\Theta^{\star q-1} \max_{1 \le j \le q} |\boldsymbol{\theta}((k+j)/n) - \boldsymbol{\theta}((k+1)/n)|. \end{aligned}$$

Let us set $\tilde\rho \in (\rho, \tau)$. Hence, for all $\boldsymbol{\theta} \in \mathcal{S}(\rho)$ and $t \in [0,1]$, we have $|\lambda|_{\max}(\Theta(t, \boldsymbol{\theta})) \le \rho$. Since $\Theta^\star < \infty$, by Lemma 12 we obtain

$$\sup_{\boldsymbol{\theta} \in \Lambda \cap \mathcal{S}} \sup_{t \in [0,1]} |\Theta^q(t, \boldsymbol{\theta})| \le C\tilde\rho^q, \qquad q \in \mathbb{N}.$$

Observe that, for all $n \ge 1$,

$$\sup_{\boldsymbol{\theta} \in \Lambda} \max_{1 \le j \le q} \max_{0 \le k \le n-q} |\boldsymbol{\theta}((k+j)/n) - \boldsymbol{\theta}((k+1)/n)| \le L(q/n)^\beta.$$

Pick $q$ and then $N$ large enough so that $C\bar{\rho}^q \leq \tau^q/2$ and $q\Theta^{\star q-1}L(q/N)^\beta \leq \tau^q/2$. The three last displays then give, for all $n \geq N$,

$$\sup_{\boldsymbol{\theta} \in \Lambda \cap \mathcal{S}} \max_{0 \leq k \leq (n-q)} |\beta_n(k+q, k; \boldsymbol{\theta})| \leq \tau^q.$$

Write $m = sq + t$, $0 \leq t < q$. For all $n \geq N, l \in \{1, \ldots, n-m\}$, and for all $\boldsymbol{\theta} \in \Lambda \cap \mathcal{S}$,

$$|\beta_n(l+m, l-1; \boldsymbol{\theta})| \leq \Theta^{\star t} \prod_{i=0}^{s-1} |\beta_n(l+(i+1)q, l+iq; \boldsymbol{\theta})|$$

$$\leq (1+\Theta^\star/\tau)^q \tau^m.$$

The proof follows. $\quad\square$

From Proposition 13 and Proposition 10, we get that under (A1),

$$(34) \qquad \sup_{(\boldsymbol{\theta}, \sigma) \in \mathcal{C}} \sup_{0 \leq k \leq n} \|\mathbf{X}_{k,n}\|_{q,\boldsymbol{\theta},\sigma} < \infty.$$

Equations (31) and (34) are referred to as uniform exponential stability and uniform $L^q$ boundedness, respectively. The bound (34) may be extended to conditional moments as follows.

PROPOSITION 14. *Assume* (A1) *with $q \geq 1$ and let $p \in [1, q]$. Let $\beta \in (0, 1]$, $L > 0$, $0 < \rho < \tau < 1$ and $0 < \sigma_- \leq \sigma_+$. Then there exists a constant $M$ such that, for all $(\boldsymbol{\theta}, \sigma) \in \mathcal{C}(\beta, L, \rho, \sigma_-, \sigma_+)$ and $1 \leq k \leq l \leq n$,*

$$(35) \qquad \mathbb{E}_{\boldsymbol{\theta},\sigma}^{\mathcal{F}_{k,n}}[|\mathbf{X}_{l,n}|^p] \leq M(1 + \tau^{(l-k)}|\mathbf{X}_{k,n}|^p),$$

*where $\mathcal{F}_{k,n} = \sigma(\mathbf{X}_{0,n}, X_{j,n}, 1 \leq j \leq k)$.*

PROOF. Equation (31) is satisfied by Proposition 10. Then under (A1) we apply Proposition 11 with $\phi(\mathbf{x}) = |\mathbf{x}|^p$, since $\phi \in \text{Li}(\mathbb{R}^d, 1, \mathbb{R}; p-1)$. Equation (35) follows from (25). $\quad\square$

*Step 2. Error decomposition.* When studying recursive algorithms of the form (3), it is convenient to rewrite the original recursion in terms of the error defined by $\delta_{k,n} := \hat{\boldsymbol{\theta}}_{k,n}(\mu) - \boldsymbol{\theta}_{k,n}$, $0 \leq k \leq n$. Let us denote, for all $\nu \geq 0$ and $\mathbf{x} \in \mathbb{R}^d$,

$$(36) \qquad \text{L}_\nu(\mathbf{x}) := \frac{\mathbf{x}}{1 + \nu|\mathbf{x}|^2} \quad \text{and} \quad \text{F}_\nu(\mathbf{x}) := \text{L}_\nu(\mathbf{x})\mathbf{x}^T.$$

The tracking error process $\{\delta_{k,n}, 0 \le k \le n\}$ obeys the following sequence of linear stochastic difference equations. For all $0 \le k < n$,

$$\delta_{k,n} = \delta_{k,n}^{(\mathrm{u})} + \delta_{k,n}^{(\mathrm{v})} + \delta_{k,n}^{(\mathrm{w})},$$

(37)
$$\delta_{k+1,n}^{(\mathrm{u})} := (I - \mu \mathrm{F}_\mu(\mathbf{X}_{k,n})) \delta_{k,n}^{(\mathrm{u})}, \qquad \delta_{0,n}^{(\mathrm{u})} = -\boldsymbol{\theta}_0,$$

$$\delta_{k+1,n}^{(\mathrm{v})} := (I - \mu \mathrm{F}_\mu(\mathbf{X}_{k,n})) \delta_{k,n}^{(\mathrm{v})} + \mu \mathrm{L}_\mu(\mathbf{X}_{k,n}) \sigma_{k+1,n} \varepsilon_{k+1,n}, \qquad \delta_{0,n}^{(\mathrm{v})} = 0,$$

$$\delta_{k+1,n}^{(\mathrm{w})} := (I - \mu \mathrm{F}_\mu(\mathbf{X}_{k,n})) \delta_{k,n}^{(\mathrm{w})} + (\boldsymbol{\theta}_{k,n} - \boldsymbol{\theta}_{k+1,n}), \qquad \delta_{0,n}^{(\mathrm{w})} = 0.$$

$\{\delta_{k,n}^{(\mathrm{u})}\}$ takes into account the way the successive estimates of the regression coefficients forget the initial error. Making a parallel with classical nonparametric function estimation, $\delta_{k,n}^{(\mathrm{w})}$ plays the role of a bias term [this term cancels when the function $t \mapsto \boldsymbol{\theta}(t)$ is constant], whereas $\delta_{k,n}^{(\mathrm{v})}$ is a stochastic disturbance. It should be stressed that the "bias term" is nondeterministic as soon as $t \mapsto \boldsymbol{\theta}(t)$ is not constant. The transient term simply writes, for all $0 \le k < n$,

(38)
$$\delta_{k+1,n}^{(\mathrm{u})} = \Psi_n(k, -1; \mu) \delta_0^{(\mathrm{u})},$$

where, for all $\mu \ge 0$ and $-1 \le j < k \le n$,

(39)
$$\Psi_n(j, j; \mu) := I \quad \text{and} \quad \Psi_n(k, j; \mu) = \prod_{l=j+1}^{k} (I - \mu \mathrm{F}_\mu(\mathbf{X}_{l,n})).$$

Let us finally define the following increment processes, for all $0 \le k < n$,

$$\xi_{k,n}^{(\mathrm{w})} := \boldsymbol{\theta}_{k,n} - \boldsymbol{\theta}_{k+1,n} \quad \text{and} \quad \xi_{k,n}^{(\mathrm{v})} := \mu \mathrm{L}_\mu(\mathbf{X}_{k,n}) \sigma_{k+1,n} \varepsilon_{k+1,n}.$$

According to these definitions, $\{\delta_{k,n}^{(\mathrm{v})}\}_{1 \le k \le n}$ and $\{\delta_{k,n}^{(\mathrm{w})}\}_{1 \le k \le n}$ obey a generic sequence of inhomogeneous stochastic recurrence equations of the form

$$\delta_{k+1,n}^{(\bullet)} = (I - \mu \mathrm{F}_\mu(\mathbf{X}_{k,n})) \delta_{k,n}^{(\bullet)} + \xi_{k,n}^{(\bullet)} = \sum_{i=0}^{k} \Psi_n(k, i; \mu) \xi_{i,n}^{(\bullet)}, \qquad 0 \le k < n.$$

(40)
In view of (38) and (40), it is clear that the stability of the product of random matrices $\Psi_n(k, i; \mu)$ plays an important role in the limiting behavior of the estimation error.

*Step* 3. *Stability of the recursive algorithm.* The following stability result for the product $\Psi_n(k, j; \mu)$ defined in (39) is an essential step for deriving risk bounds for the estimator $\hat{\boldsymbol{\theta}}_n$.

THEOREM 15. *Assume* (A1) *with* $q \geq 4$. *Let* $\beta \in (0, 1]$, $L > 0$, $0 < \rho < 1$, *and* $0 < \sigma_- \leq \sigma_+$. *Then for all* $p \geq 1$ *there exist constants* $M, \delta > 0$ *and* $\mu_0 > 0$, *such that, for all* $0 \leq j \leq k \leq n$, $\mu \in [0, \mu_0]$ *and* $(\boldsymbol{\theta}, \sigma) \in \mathcal{C}(\beta, L, \rho, \sigma_-, \sigma_+)$,

$$(41) \qquad \|\Psi_n(k, j; \mu)\|_{p, \boldsymbol{\theta}, \sigma} \leq M(1 - \delta\mu)^{k-j}.$$

Similar stability results have been obtained in the framework of classical recursive estimation algorithms (see, e.g., [13, 20]), but cannot be applied directly to our nonstationary and nonparametric context. Let us sketch the main arguments of the proof in a more general context. Let $\{A_k(\nu), k \geq 0, \nu > 0\}$ be an $\mathbb{M}_d^+$-valued process such that

(C-1) for all $k \in \mathbb{N}$ and $\nu \in [0, \nu_1]$, $|\nu A_k(\nu)| \leq 1$.

Here $A_0(\nu), A_1(\nu), \ldots$ correspond to the matrices $\mathrm{F}_\nu(\mathbf{X}_{j+1,n}), \mathrm{F}_\nu(\mathbf{X}_{j+2,n}), \ldots$, which appear in the product $\Psi_n(k, j; \nu)$ for some fixed $j$, $k$ and $n$. Taking $p = 1$ in (41) [this can be done without loss of generality under (C-1)], we want to prove that $\mathbb{E}|\prod_{k=1}^l (I - \nu A_k(\nu))| \leq C(1 - \delta\nu)^l$ for some positive constants $C$ and $\delta$.

First observe that it is sufficient to have $\mathbb{E}^{\mathcal{F}_k}[|I - \nu A_{k+1}(\nu)|] \leq 1 - \delta\nu$ or, more generally, for some fixed integer $r$ and $\alpha > 0$,

$$(42) \qquad \mathbb{E}^{\mathcal{F}_k}\left[\left|\prod_{l=k+1}^{k+r}(I - \nu A_l(\nu))\right|\right] \leq 1 - \alpha\nu.$$

To obtain this inequality, we expand the product and then use Lemma C.1 so that

$$\left|\prod_{l=k+1}^{k+r}(I - \nu A_l(\nu))\right| \leq 1 - \nu\lambda_{\min}\left(\sum_{l=k+1}^{k+r} A_l(\nu)\right) + \mathcal{R},$$

where $\mathcal{R}$ is some remainder term which will be controlled by using $\sum_{l=k+1}^{k+r} |A_l(\nu)|^s$ for some $s > 1$. It turns out that, in the previous display, the conditional expectation of the RHS given $\mathcal{F}_k$ may be controlled only on a set $\{\phi_k \leq R_1\}$, where $R_1 > 0$ and $\{\phi_k\}$ is a positive adapted sequence (in the TVAR context, $\phi_k$ corresponds to $|\mathbf{X}_{k,n}|$). This yields the following two conditions:

(C-2) there exists $\alpha_1 > 0$ such that, for all $k \in \mathbb{N}$ and $\nu \in [0, \nu_1]$,

$$\mathbb{E}^{\mathcal{F}_k}\left[\lambda_{\min}\left(\sum_{l=k+1}^{k+r} A_l(\nu)\right)\right] \geq \alpha_1 \mathbf{I}(\phi_k \leq R_1);$$

(C-3) there exist $s > 1$ and $C_1 > 0$ such that, for all $k \in \mathbb{N}$ and $\nu \in [0, \nu_1]$,

$$\mathbf{I}(\phi_k \leq R_1) \sum_{l=k+1}^{k+r} \mathbb{E}^{\mathcal{F}_k}[|A_l(\nu)|^s] \leq C_1.$$

In turn, (42) may be shown only on the set $\{\phi_k \leq R_1\}$, and it remains to check that this happens for sufficiently many $k$'s. For this we use a classical Lyapounov condition, namely,

(C-4) there exist $\lambda < 1$, $B \geq 1$ and an adapted process $\{V_k, k \geq 0\}$ on $[1, \infty)$ such that, for all $k \in \mathbb{N}$,

$$\mathbb{E}^{\mathcal{F}_k}[V_{k+r}] \leq \lambda V_k \mathbf{I}(\phi_k > R_1) + BV_k \mathbf{I}(\phi_k \leq R_1).$$

The previous arguments yield the following general result, whose precise technical proof is postponed to Appendix A for convenience.

THEOREM 16. Let $(\Omega, \mathcal{F}, \mathbb{P}, \{\mathcal{F}_k : k \in \mathbb{N}\})$ be a filtered space. Let $\{\phi_k, k \geq 0\}$ be a nonnegative adapted process and, for any $\nu \geq 0$, let $A(\nu) := \{A_k(\nu), k \geq 0\}$ be an adapted $\mathbb{M}_d^+$-valued process. Let $r \geq 1$, $R_1 > 0$ and $\nu_1 > 0$ such that (C-1)–(C-4) hold. Then, for any $p \geq 1$ there exist $C_0 > 0$, $\delta_0 > 0$ and $\nu_0 > 0$ depending only on $p, \nu_1, \alpha_1, R_1, C_1, s, B$ and $\lambda$ such that, for all $\nu \in [0, \nu_0]$ and $n \geq 1$,

$$\mathbb{E}^{\mathcal{F}_0}\left[\left|\prod_{i=1}^n (I - \nu A_i(\nu))\right|^p\right] \leq C_0 e^{-\delta_0 \nu n} V_0.$$

Having this general result in hand, we now prove Theorem 15.

PROOF OF THEOREM 15. To verify (C-2), which is referred to as the *persistence of excitation* property in the control theory literature, we need some intermediary results which hold under the assumptions of Theorem 15.

LEMMA 17. There exists $C > 0$ such that, for all $d \leq j \leq n$ and $(\boldsymbol{\theta}, \sigma) \in \mathcal{C}(\beta, L, \rho, \sigma_-, \sigma_+)$,

(43) $$\lambda_{\min}(\mathbb{E}_{\boldsymbol{\theta}, \sigma}[\mathrm{F}_\nu(\mathbf{X}_{j,n})]) \geq C, \qquad \nu \in [0, 1],$$

and, for all $0 \leq k \leq j - d$,

(44) $$\lambda_{\min}(\mathbb{E}_{\boldsymbol{\theta}, \sigma}^{\mathcal{F}_{k,n}}[\mathrm{F}_0(\mathbf{X}_{j,n})]) \geq C,$$

(45) $$\lambda_{\min}(\mathbb{E}_{\boldsymbol{\theta}, \sigma}^{\mathcal{F}_{k,n}}[\mathrm{F}_\nu(\mathbf{X}_{j,n})]) \geq \frac{C}{1 + |\mathbf{X}_{k,n}|^4}, \qquad \nu \in [0, 1].$$

PROOF. Let $\mathbf{x} \in \mathbb{R}^d$, $|\mathbf{x}| = 1$, $d \leq j \leq n$, and $(\boldsymbol{\theta}, \sigma) \in \mathcal{C} := \mathcal{C}(\beta, L, \rho, \sigma_-, \sigma_+)$. Write

$$\mathbf{X}_{j,n} = \beta_n(j, j-d; \boldsymbol{\theta})\mathbf{X}_{j-d,n} + \sum_{i=1}^d \beta_n(j, j-d+i; \boldsymbol{\theta})\boldsymbol{\sigma}_{j-d+i,n}\varepsilon_{j-d+i,n},$$

where $\beta_n$ is defined by (33). We have $\mathbb{E}_{\boldsymbol{\theta},\sigma}^{\mathcal{F}_{j-d,n}}[(\mathbf{x}^T\mathbf{X}_{j,n})^2] \geq \sigma^2\mathbf{x}^T H_{j,n,d}(\boldsymbol{\theta})\mathbf{x}$, where $H_{j,n,d}(\boldsymbol{\theta}) := C_{j,n,d}(\boldsymbol{\theta})C_{j,n,d}^T(\boldsymbol{\theta})$ is the controllability Gramian (see [17])

$$C_{j,n,d}(\boldsymbol{\theta}) := [\beta_n(j,j-d+1;\boldsymbol{\theta})G \cdots \beta_n(j,j;\boldsymbol{\theta})G], \qquad \text{where } G := [1\,0\ldots 0]^T.$$

One easily shows that, for $i = 1,\ldots,d$, $\beta_n(j,j+1-i;\boldsymbol{\theta})G$ has a unit $i$th coordinate and zero coordinates below. Hence, $\det(H_{j,n,d}(\theta)) = \det(C_{j,n,d}(\theta)) = 1$. In addition, from exponential stability we have, for all $(\boldsymbol{\theta},\sigma) \in \mathcal{C}$, $|C_{j,n,d}(\boldsymbol{\theta})| \leq M$ for some positive $M$ not depending on $(j,n)$. Hence, for all $d \leq j \leq n$ and $(\boldsymbol{\theta},\sigma) \in \mathcal{C}$,

$$\lambda_{\min}(H_{j,n,d}(\boldsymbol{\theta})) \geq \frac{\det(H_{j,n,d}(\boldsymbol{\theta}))}{\lambda_{\max}^{d-1}(H_{j,n,d}(\boldsymbol{\theta}))} \geq M^{-(d-1)}.$$

It follows that, for all $0 \leq k \leq j-d$, $\mathbf{x} \in \mathbb{R}^d$ such that $|\mathbf{x}| = 1$ and $(\boldsymbol{\theta},\sigma) \in \mathcal{C}$,

$$(46) \qquad \mathbf{x}^T\mathbb{E}_{\boldsymbol{\theta},\sigma}^{\mathcal{F}_{k,n}}[F_0(\mathbf{X}_{j,n})]\mathbf{x} = \mathbb{E}_{\boldsymbol{\theta},\sigma}^{\mathcal{F}_{k,n}}[\mathbb{E}_{\boldsymbol{\theta},\sigma}^{\mathcal{F}_{j-d,n}}[(\mathbf{x}^T\mathbf{X}_{j,n})^2]] \geq \sigma_-^2 M^{-(d-1)},$$

showing (44) for any $C \geq \sigma_-^2 M^{-(d-1)}$. Equation (43) also follows for $\nu = 0$. For all $0 \leq k \leq j-d$ and for all $\nu \in [0,1]$, write

$$\mathbb{E}_{\boldsymbol{\theta},\sigma}^{\mathcal{F}_{k,n}}[(\mathbf{x}^T\mathbf{X}_{j,n})^2] = \mathbb{E}_{\boldsymbol{\theta},\sigma}^{\mathcal{F}_{k,n}}\left[\frac{|\mathbf{x}^T\mathbf{X}_{j,n}|}{(1+\nu|\mathbf{X}_{j,n}|^2)^{1/2}}\{|\mathbf{x}^T\mathbf{X}_{j,n}|(1+\nu|\mathbf{X}_{j,n}|^2)^{1/2}\}\right].$$

The Cauchy–Schwarz inequality shows that

$$\mathbf{x}^T\mathbb{E}_{\boldsymbol{\theta},\sigma}^{\mathcal{F}_{k,n}}[F_\nu(\mathbf{X}_{j,n})]\mathbf{x} \geq \frac{(\mathbb{E}_{\boldsymbol{\theta},\sigma}^{\mathcal{F}_{k,n}}[(\mathbf{x}^T\mathbf{X}_{j,n})^2])^2}{\mathbb{E}_{\boldsymbol{\theta},\sigma}^{\mathcal{F}_{k,n}}[(\mathbf{x}^T\mathbf{X}_{j,n})^2(1+\nu|\mathbf{X}_{j,n}|^2)]}.$$

Since $|\mathbf{x}^T\mathbf{X}_{j,n}| \leq |\mathbf{X}_{j,n}|$ for $|\mathbf{x}| = 1$, by applying (46) we get, for all $\nu \in [0,1]$,

$$\lambda_{\min}(\mathbb{E}_{\boldsymbol{\theta},\sigma}^{\mathcal{F}_{k,n}}[F_\nu(\mathbf{X}_{j,n})]) \geq \frac{(\sigma_-^2 M^{-(d-1)})^2}{\mathbb{E}_{\boldsymbol{\theta},\sigma}^{\mathcal{F}_{k,n}}[|\mathbf{X}_{j,n}|^2] + \mathbb{E}_{\boldsymbol{\theta},\sigma}^{\mathcal{F}_{k,n}}[|\mathbf{X}_{j,n}|^4]}.$$

The proof of (45) then follows from (35). The proof of (43) is along the same lines. □

LEMMA 18. *Let* $\phi \in \mathrm{Li}(\mathbb{R}^d, 1, \mathbb{R}; 1)$. *There exists* $M > 0$ *such that, for all* $1 \leq i \leq k \leq n$, $\nu \geq 0$ *and* $(\boldsymbol{\theta},\sigma) \in \mathcal{C}(\beta, L, \rho, \sigma_-, \sigma_+)$,

$$\mathbb{E}_{\boldsymbol{\theta},\sigma}^{\mathcal{F}_{i,n}}\left[\left|\sum_{j=i+1}^{k}(\phi(\mathbf{X}_{j,n}) - \mathbb{E}_{\boldsymbol{\theta},\sigma}^{\mathcal{F}_{i,n}}[\phi(\mathbf{X}_{j,n})])\right|^2\right] \leq M(k-i)|\phi|_{\mathrm{Li}(1)}^2(1+|\mathbf{X}_{i,n}|^4).$$

PROOF. For $j \in \{i+1,\ldots,k\}$, denote $\Delta_j = \phi(\mathbf{X}_{j,n}) - \mathbb{E}_{\boldsymbol{\theta},\sigma}^{\mathcal{F}_{i,n}}[\phi(\mathbf{X}_{j,n})]$. For all $i \leq j \leq l \leq n$, $\mathbb{E}_{\boldsymbol{\theta},\sigma}^{\mathcal{F}_{j,n}}[\Delta_l] = \mathbb{E}_{\boldsymbol{\theta},\sigma}^{\mathcal{F}_{j,n}}[\phi(\mathbf{X}_{l,n})] - \mathbb{E}_{\boldsymbol{\theta},\sigma}^{\mathcal{F}_{i,n}}[\phi(\mathbf{X}_{l,n})]$. The model is uniformly exponentially and $L^q$ stable [see Proposition 13 and (34)]. From Proposition 11 there exist $\tau \in (\rho, 1)$ and $C_1, C_2 > 0$ such that, for all $0 \le i \le j \le l \le n$, $(\boldsymbol{\theta}, \sigma) \in \mathcal{C}(\beta, L, \rho, \sigma_-, \sigma_+)$,

$$|\mathbb{E}_{\boldsymbol{\theta}, \sigma}^{\mathcal{F}_{j,n}}[\Delta_l]| \le C_1 \tau^{l-j} |\phi|_{\mathrm{Li}(1)} (1 + |\mathbf{X}_{j,n}|^2 + |\mathbf{X}_{i,n}|^2),$$

$$|\mathbb{E}_{\boldsymbol{\theta}, \sigma}^{\mathcal{F}_{j,n}}[\Delta_l^2]| \le C_2 |\phi|_{\mathrm{Li}(1)}^2 (1 + |\mathbf{X}_{i,n}|^4).$$

From Proposition 14 we get that there exists $C_3 > 0$ such that, for all $0 \le i \le j \le l \le n$, $(\boldsymbol{\theta}, \sigma) \in \mathcal{C}(\beta, L, \rho, \sigma_-, \sigma_+)$,

$$|\mathbb{E}_{\boldsymbol{\theta}, \sigma}^{\mathcal{F}_{i,n}}[\Delta_j \Delta_l]| = |\mathbb{E}_{\boldsymbol{\theta}, \sigma}^{\mathcal{F}_{i,n}}[\Delta_j \mathbb{E}_{\boldsymbol{\theta}, \sigma}^{\mathcal{F}_{j,n}}[\Delta_l]]|$$

$$\le (\mathbb{E}_{\boldsymbol{\theta}, \sigma}^{\mathcal{F}_{i,n}}[\Delta_j^2])^{1/2} (\mathbb{E}_{\boldsymbol{\theta}, \sigma}^{\mathcal{F}_{i,n}}[(\mathbb{E}_{\boldsymbol{\theta}, \sigma}^{\mathcal{F}_{j,n}}[\Delta_l])^2])^{1/2}$$

$$\le C_3 |\phi|_{\mathrm{Li}(1)}^2 \tau^{l-j} (1 + |\mathbf{X}_{i,n}|^4),$$

and the result follows. $\square$

LEMMA 19. *For all $R_1, \alpha_1 > 0$, there exists $r_0 \ge 1$ such that, for all $(\boldsymbol{\theta}, \sigma) \in \mathcal{C}(\beta, L, \rho, \sigma_-, \sigma_+)$, $r \ge r_0$, $n \ge r$, $k = 0, \ldots, n-r$ and $\nu \in (0, 1]$,*

$$(47) \qquad \mathbb{E}_{\boldsymbol{\theta}, \sigma}^{\mathcal{F}_{k,n}} \left[ \lambda_{\min} \left( \sum_{l=k+1}^{k+r} \mathrm{F}_\nu(\mathbf{X}_{l,n}) \right) \right] \ge \alpha_1 \mathbf{I}(|\mathbf{X}_{k,n}| \le R_1),$$

*where $\mathrm{F}_\nu$ is defined in (36). In addition, there exist constants $\delta > 0$ and $\mu_0 > 0$, such that, for all $(\boldsymbol{\theta}, \sigma) \in \mathcal{C}(\beta, L, \rho, \sigma_-, \sigma_+)$, $d \le k \le n$ and $\nu \in [0, \nu_0]$,*

$$(48) \qquad |I - \nu \mathbb{E}_{\boldsymbol{\theta}, \sigma}[\mathrm{F}_\nu(\mathbf{X}_{k,n})]| \le 1 - \delta \nu.$$

PROOF. For any symmetric matrix $A$, we have $|\lambda_{\min}(A)| \le |\lambda|_{\max}(A) = |A|$ (recall that $|\cdot|$ denotes the operator norm) and $\lambda_{\min}(A) = \inf_{|\mathbf{x}|=1} x^T A x$. From the last assertion, it follows that, for any symmetric matrix $B$ having the same size as $A$, $\lambda_{\min}(A + B) \ge \lambda_{\min}(A) + \lambda_{\min}(B)$. Therefore,

$$\lambda_{\min}(A) \ge \lambda_{\min}(B) + \lambda_{\min}(A - B) \ge \lambda_{\min}(B) - |A - B|.$$

Applying these elementary facts, we get, for all $0 \le k < k + r \le n$,

$$\lambda_{\min} \left( \sum_{j=k+1}^{k+r} \mathrm{F}_\nu(\mathbf{X}_{j,n}) \right)$$

$$\ge \sum_{j=k+1}^{k+r} \lambda_{\min}(\mathbb{E}_{\boldsymbol{\theta}, \sigma}^{\mathcal{F}_{k,n}}[\mathrm{F}_\nu(\mathbf{X}_{j,n})]) - \left| \sum_{j=k+1}^{k+r} (\mathrm{F}_\nu(\mathbf{X}_{j,n}) - \mathbb{E}_{\boldsymbol{\theta}, \sigma}^{\mathcal{F}_{k,n}}[\mathrm{F}_\nu(\mathbf{X}_{j,n})]) \right|.$$

From its definition in (36), $\mathrm{F}_\nu(\mathbf{x}) \in \mathbb{M}_d^+$ for all $\mathbf{x} \in \mathbb{R}^d$ and $\nu \ge 0$, and

$$(49) \qquad \sup_{\nu \ge 0} |\mathrm{F}_\nu|_{\mathrm{Li}(1)} < \infty.$$

Applying Lemmas 17 and 18, we obtain that there exist $C, M > 0$ such that, for all $(\boldsymbol{\theta}, \sigma) \in \mathcal{C}$ and $0 \leq k < k + r \leq n$,

$$\mathbb{E}_{\boldsymbol{\theta}, \sigma}^{\mathcal{F}_{k,n}}\left[\lambda_{\min}\left(\sum_{j=k+1}^{k+r} \mathrm{F}_\nu(\mathbf{X}_{j,n})\right)\right] \geq \frac{C(r-d)}{1+|\mathbf{X}_{k,n}|^4} - M\sqrt{r}(1+|\mathbf{X}_{k,n}|^4)^{1/2}.$$

Now pick two positive numbers $R_1$ and $\alpha_1$. If $|\mathbf{X}_{k,n}| > R_1$, equation (47) is clearly satisfied. If $|\mathbf{X}_{k,n}| \leq R_1$, the last equation implies that, for all $0 \leq k < k + r \leq n$,

$$\mathbb{E}_{\boldsymbol{\theta}, \sigma}^{\mathcal{F}_{k,n}}\left[\lambda_{\min}\left(\sum_{j=k+1}^{k+r} \mathrm{F}_\nu(\mathbf{X}_{j,n})\right)\right] \geq \frac{C(r-d)}{1+R_1} - M\sqrt{r}(1+R_1^4)^{1/2}.$$

We may thus find $r_0$ such that the RHS of this inequality is greater than or equal to $\alpha_1$ for all $r \geq r_0$. This concludes the proof of (47).

From the uniform $L^2$-boundedness and (49), we get that there exists $M$ such that, for all $(\boldsymbol{\theta}, \sigma) \in \mathcal{C}$, $\nu \in [0,1]$, $d \leq k \leq n$, $|\mathbb{E}_{\boldsymbol{\theta}, \sigma}[\mathrm{F}_\nu(\mathbf{X}_{k,n})]| \leq M$. Thus, using Lemma C.1, for all $\nu \in [0, 1/M]$, $d \leq k \leq n$, and $(\boldsymbol{\theta}, \sigma) \in \mathcal{C}$,

$$|I - \nu\mathbb{E}_{\boldsymbol{\theta}, \sigma}[\mathrm{F}_\nu(\mathbf{X}_{k,n})]| = 1 - \nu\lambda_{\min}(\mathbb{E}_{\boldsymbol{\theta}, \sigma}[\mathrm{F}_\nu(\mathbf{X}_{k,n})])$$

and the proof of (48) follows from Lemma 17. $\square$

We now turn back to the proof of Theorem 15 by applying Theorem 16 to the sequence

$$\{(A_l = \mathrm{F}_\nu(\mathbf{X}_{j+1+l,n}), \mathcal{F}_{j+1+l,n}), l = 0, \ldots, k-j-1\}.$$

It thus remains to show that conditions (C-1)–(C-4) of Theorem 16 hold with constants $r, R_1, \nu_1, \alpha_1, C_1, \lambda, B$ and $s$, which neither depend on $(\boldsymbol{\theta}, \sigma) \in \mathcal{C}$ nor on $j, k, n$. Set $\alpha_1 = 1$ and, for $l = 0, \ldots, k-j-1$, $V_l := 1 + |\mathbf{X}_{j+l,n}|$ and $\phi_l := |\mathbf{X}_{j+l,n}|$.

Condition (C-1). For all $\nu \geq 0$ and $\mathbf{x} \in \mathbb{R}^d$, $\nu|\mathrm{F}_\nu(\mathbf{x})| = \nu|\mathbf{x}|/(1+\nu|\mathbf{x}|) \leq 1$, which yields (C-1).

Condition (C-2). From Lemma 19, we may choose $r_0$ depending only on $R_1$ and $\alpha_1$ such that (C-2) holds for all $r \geq r_0$.

Condition (C-3). From Lemma C.2, $\sup_{\nu \geq 0}||\mathrm{F}_\nu|^{q/2}|_{\mathrm{Li}(q-1)} < \infty$. From Proposition 14, there exists $M > 0$ such that, for all $\nu \geq 0$, $0 \leq i < i + l \leq n$ and $(\boldsymbol{\theta}, \sigma) \in \mathcal{C}$,

$$\mathbb{E}_{\boldsymbol{\theta}, \sigma}^{\mathcal{F}_{i,n}}[|\mathrm{F}_\nu(\mathbf{X}_{i+l,n})|^{q/2}] \leq M(1 + \tau^{ql}|\mathbf{X}_{i,n}|^q).$$

Hence, (C-3) is obtained with $s = q/2$ and $C_1 = Mr(1 + R_1^q)$.

Condition (C-4). Let $\tau \in (\rho, 1)$. From Proposition 14, there exists $M$ such that, for all $0 \le i < i+r \le n$ and $(\boldsymbol{\theta}, \sigma) \in \mathcal{C}$, $\mathbb{E}_{\boldsymbol{\theta}, \sigma}^{\mathcal{F}_{i,n}}[|\mathbf{X}_{i+r,n}|] \le M(1 + \tau^r |\mathbf{X}_{i,n}|)$; thus, for any $R_1 > 0$,

$$\mathbb{E}_{\boldsymbol{\theta}, \sigma}^{\mathcal{F}_{i,n}}[1 + |\mathbf{X}_{i+r,n}|] \le M\tau^r |\mathbf{X}_{i,n}| + M + 1$$
$$\le \left( M\tau^r + \frac{M+1}{R_1} \mathbf{I}(|\mathbf{X}_{i,n}| > R_1) \right)$$
$$\times (1 + |\mathbf{X}_{i,n}|) + (M+1)\mathbf{I}(|\mathbf{X}_{i,n}| \le R_1).$$

Choose $r \ge r_0$ and $R_1 > 0$ so that $M\tau^r + (M+1)/R_1 < 1$. Condition (C-4) is then satisfied with $\lambda := M\tau^r + (M+1)/R_1 < 1$ and $B = 1 + M$.

Finally, we obtain, for some positive constants $r, C_0, \delta_0$ and $\mu_0$, for all $\nu \in (0, \mu_0]$, $(\boldsymbol{\theta}, \sigma) \in \mathcal{C}$ and $0 \le j < k \le n$ such that $n - j \ge r$,

$$\|\Psi_n(k, j; \nu)\|_{p, \boldsymbol{\theta}, \sigma}^p = \mathbb{E}[\mathbb{E}_{\boldsymbol{\theta}, \sigma}^{\mathcal{F}_{j+1,n}}[\Psi_n(k, j; \nu)]]$$
$$\le C_0 e^{-\delta_0 \nu n}(1 + \|\mathbf{X}_{j+1,n}\|_{p, \boldsymbol{\theta}, \sigma}).$$

The uniform boundedness (34) then yields (41) when $n - j \ge r$. The restriction $n - j \ge r$ above is needed because (C-4), (C-2) and (C-3) are well defined only for $n - j < r$. Now recall that (C-1) implies $|\Psi_n(k, j; \nu)| \le 1$. The result for $n - j < r$ (implying $k - j < r$) follows by taking $M \ge (1 - \delta\nu_1)^{-r}$ in (41). $\square$

*Step* 4. *Error bounds.* Similarly to (49), one easily shows that

$$(50) \qquad \sup_{\nu \ge 0} |\mathrm{L}_\nu|_{\mathrm{Li}(0)} < \infty,$$

where L is defined in (36). The $L^q$-boundedness of $\{\mathbf{X}_{k,n}, 1 \le k \le n\}$ then gives

$$(51) \qquad \mathbf{F}_{q/2}^\star := \sup_{(\boldsymbol{\theta}, \sigma) \in \mathcal{C}} \sup_{\nu \ge 0} \sup_{0 \le k \le n} \|\mathrm{F}_\nu(\mathbf{X}_{k,n})\|_{q/2, \boldsymbol{\theta}, \sigma} < \infty,$$

$$(52) \qquad \mathrm{L}_q^\star := \sup_{(\boldsymbol{\theta}, \sigma) \in \mathcal{C}} \sup_{\nu \ge 0} \sup_{0 \le k \le n} \|\mathrm{L}_\nu(\mathbf{X}_{k,n})\|_{q, \boldsymbol{\theta}, \sigma} < \infty.$$

From now on, for convenience, we let the same $M$, $\delta$ and $\mu_0$ denote positive constants depending neither on the time indices $i, j, k, n, \ldots$, the step-size $\mu$ nor on the parameter $(\boldsymbol{\theta}, \sigma)$.

Applying (38) and (41), for all $u \ge 1$, there exist $M > 0$ and $\mu_0 > 0$ such that, for all $\mu \in (0, \mu_0]$, $(\boldsymbol{\theta}, \sigma) \in \mathcal{C}$ and $1 \le k \le n$,

$$(53) \qquad \|\delta_{k,n}^{(\mathrm{u})}\|_{u, \boldsymbol{\theta}, \sigma} \le M(1 - \delta\mu)^k |\boldsymbol{\theta}(0)|.$$

Define

$$\Xi_n^{(\bullet)}(k,k) := 0 \quad \text{and} \quad \Xi_n^{(\bullet)}(k,j) := \sum_{i=j}^{k-1} \xi_{i,n}^{(\bullet)}, \qquad 0 \le j < k \le n.$$

For all $1 \le j \le k \le n$, we have

$$\Psi_n(k-1,j;\mu) - \Psi_n(k-1,j-1;\mu) = \mu \Psi_n(k-1,j;\mu) F_\mu(\mathbf{X}_{j-1,n}).$$

By integration by parts, for all $1 \le k \le n$ and $\mu \in (0,\mu_0]$, (40) reads

$$\begin{aligned}
(54) \qquad \delta_{k,n}^{(\bullet)} &= \Psi_n(k-1,0;\mu) \Xi_n^{(\bullet)}(k,0) \\
&\quad + \mu \sum_{j=1}^{k-1} \Psi_n(k-1,j;\mu) F_\mu(\mathbf{X}_{j-1,n}) \Xi_n^{(\bullet)}(k,j).
\end{aligned}$$

By applying the Hölder inequality and using (41) and (51), we get that, for any $u \in (1, q/2)$, there exists $M > 0$ such that, for all $\mu \in [0, \mu_0]$, $(\boldsymbol{\theta}, \sigma) \in \mathcal{C}$ and $0 \le j \le k \le n$,

$$(55) \qquad \|\Psi_n(k,j;\mu) F_\mu(\mathbf{X}_{j-1,n})\|_{u,\boldsymbol{\theta},\sigma} \le M(1-\delta\mu)^{k-j+1}.$$

We consider now the two terms $\delta_{k,n}^{(\mathrm{w})}$ and $\delta_{k,n}^{(\mathrm{v})}$ separately. We have, for all $0 \le j < k \le n$,

$$(56) \qquad |\Xi_n^{(\mathrm{w})}(k,j)| := \left| \sum_{i=j}^{k-1} \xi_{i,n}^{(\mathrm{w})} \right| = |\boldsymbol{\theta}_{j,n} - \boldsymbol{\theta}_{k,n}| \le |\boldsymbol{\theta}|_{\Lambda,\beta} n^{-\beta} (k-j)^\beta.$$

Inserting (56), (41) and (55) into (54), we thus get that there exist $M > 0$, $\delta > 0$ and $\mu_0 > 0$ such that, for all $(\boldsymbol{\theta}, \sigma) \in \mathcal{C}$, $1 \le k \le n$ and $\mu \in (0, \mu_0]$,

$$\|\delta_{k,n}^{(\mathrm{w})}\|_{u,\boldsymbol{\theta},\sigma} \le M\left( (1-\delta\mu)^k (k/n)^\beta + \mu \sum_{j=1}^{k-1} (1-\delta\mu)^{k-j} ((k-j)/n)^\beta \right).$$

By Lemma C.3, picking $u \ge p$, we obtain, for all $\mu \in (0, \mu_0]$ and $(\boldsymbol{\theta}, \sigma) \in \mathcal{C}$,

$$(57) \qquad \|\delta_{k,n}^{(\mathrm{w})}\|_{p,\boldsymbol{\theta},\sigma} \le M |\boldsymbol{\theta}|_{\Lambda,\beta} (\mu n)^{-\beta}, \qquad 1 \le k \le n.$$

We finally bound $\delta^{(\mathrm{v})}$. Note that, for each $n > 1$, $\{\sigma_{i,n} L_\mu(\mathbf{X}_{i,n})\varepsilon_{i+1,n}, i = 1, \ldots, n-1\}$ is an $\mathcal{F}_{i+1,n}$-adapted martingale increment sequence. The Burkholder inequality (see [14], Theorem 2.12) and (52) give, for all $(\boldsymbol{\theta}, \sigma) \in \mathcal{C}$ and $\mu \ge 0$,

$$(58) \qquad \left\| \sum_{i=j}^k \sigma_{i,n} L_\mu(\mathbf{X}_{i,n}) \varepsilon_{i+1,n} \right\|_{q,\boldsymbol{\theta},\sigma} \le M(k-j+1)^{1/2}, \qquad 1 \le j \le k < n.$$

By Lemma C.3, using (54), (58) and (55) with $u$ such that $1/p = 1/u + 1/q$, we get, for all $\mu \in (0, \mu_0]$ and $(\boldsymbol{\theta}, \sigma) \in \mathcal{C}$,

$$(59) \qquad \|\delta_{k,n}^{(\mathrm{v})}\|_{p,\boldsymbol{\theta},\sigma} \leq M\sqrt{\mu}, \qquad 1 \leq k \leq n.$$

Equation (13) easily follows from (37), (53), (57) and (59) and Theorem 2 is obtained.

**5. Proof of Theorem 4.** By writing $\hat{\delta}_n := \mathbf{u}^T \hat{\boldsymbol{\delta}}_n$, (14) simply means that, for all real-valued estimators $\hat{\delta}_n(\mathbf{X}_{0,n}, X_{1,n}, \ldots, X_{n,n})$ and $\mathbf{u} = (u_1, \ldots, u_d) \in \mathbb{R}^d$ such that $|\mathbf{u}| = 1$,

$$\sup_{(\boldsymbol{\theta}, \sigma) \in \mathcal{C}} \mathbb{E}_{\boldsymbol{\theta}, \sigma}[|\hat{\delta}_n - \mathbf{u}^T \boldsymbol{\theta}(t)|^2] \geq \alpha n^{-2\beta/(1+2\beta)}.$$

Denote by $\lfloor \beta \rfloor$ the largest integer strictly smaller than $\beta$. Let $\phi : \mathbb{R} \to \mathbb{R}$ be a $C^\infty$ symmetric function decreasing on $\mathbb{R}_+$ such that $\phi(0) = 1$, $\phi(u) = 0$ for all $|u| \geq 1$ and

$$(60) \qquad \sup_{x,y \in \mathbb{R}, x \neq y} \frac{|\phi^{(\lfloor \beta \rfloor)}(x) - \phi^{(\lfloor \beta \rfloor)}(y)|}{|x-y|^{\beta - \lfloor \beta \rfloor}} \leq 1.$$

Let $\lambda : \mathbb{R} \to \mathbb{R}_+$ be a $C^1$ p.d.f. with respect to the Lebesgue measure vanishing outside the segment $[-1, +1]$ and such that

$$(61) \qquad \int_{-1}^{1} \left( \frac{\dot{\lambda}(x)}{\lambda(x)} \right)^2 \lambda(x) \, dx < \infty,$$

where $\dot{\lambda}$ is the derivative of the density $\lambda$. Let $(v_n)$ and $(w_n)$ be two nondecreasing sequences of positive numbers to be specified later such that

$$(62) \qquad \lim_{n \to \infty} (v_n^{-1} + w_n^{-1} + n^{-1} w_n) = 0 \quad \text{and} \quad \sup_{n \geq 0} v_n^{-1} w_n^\beta \leq L.$$

Let $t \in [0, 1]$ and $\mathbf{u} \in \mathbb{R}^d$ such that $|\mathbf{u}| = 1$. Define $\phi_{n,t} : [0, 1] \to \mathbb{R}^d$, $s \mapsto \phi_{n,t}(s) := \phi((s-t)w_n)\mathbf{u}$ and let $\sigma : [0, 1] \to \mathbb{R}^+$, $s \mapsto \sigma(s) = \sigma_+$. For $n \geq 1$, define

$$(63) \qquad X_{k+1,n} = \eta\phi((k/n - t)w_n)\mathbf{u}^T \mathbf{X}_{k,n} + (\sigma_+)\varepsilon_{k+1,n}, \qquad 0 \leq k \leq n-1.$$

From (A2), it follows that, for all $0 \leq k \leq n-1$,

$$x \mapsto p_{k+1,n}((x - \eta\phi_{n,t}(k/n)\mathbf{u}^T \mathbf{y})/\sigma_+)/\sigma_+$$

is the conditional density of $X_{k+1,n}$ given $\mathbf{X}_{k,n} = \mathbf{y}$ and parameter $\eta$. Since the distribution of $\mathbf{X}_{0,n}$ does not depend on $\eta$, and using the conditional densities above to compute the joint density of $\{\mathbf{X}_{0,n}, X_{1,n}, \ldots, X_{n,n}\}$, the

Fisher information associated with the one-dimensional parametric model defined by (63) is

$$\mathcal{I}_n(\eta) = \sigma_+^{-2} \mathbb{E}_{\eta\phi_{n,t},\sigma}\left[\left(\sum_{k=1}^n \phi_{n,t}(k/n)\mathbf{u}^T\mathbf{X}_{k-1,n}\frac{\dot{p}_{k,n}}{p_{k,n}}(\varepsilon_{k,n})\right)^2\right].$$

Now, under (A2), the summand in this equation is a martingale increments sequence whose variances are bounded by $\phi_{n,t}^2(k/n)\mathbb{E}_{\eta\phi_{n,t},\sigma}[(u^T\mathbf{X}_{k-1,n})^2]\mathcal{I}_\varepsilon$. Since $|u| = 1$, $(u^T\mathbf{X})^2 \le |\mathbf{X}|^2$, and we finally obtain

$$(64) \qquad \mathcal{I}_n(\eta) \le \sigma_+^{-2}\mathcal{I}_\varepsilon \sum_{k=1}^n \phi_{n,t}^2(k/n)\mathbb{E}_{\eta\phi_{n,t},\sigma}[|\mathbf{X}_{k-1,n}|^2].$$

From (60) and (62), for all $\eta \in [-v_n^{-1}, v_n^{-1}]$,

$$|\eta\phi_{n,t}|_{\Lambda,\beta} \le \frac{w_n^{\lfloor\beta\rfloor}}{v_n}\sup_{0\le s<s'\le 1}\frac{|\phi^{(\lfloor\beta\rfloor)}((s'-t)w_n) - \phi^{(\lfloor\beta\rfloor)}((s-t)w_n)|}{|s'-s|^{\beta-\lfloor\beta\rfloor}} \le \frac{w_n^\beta}{v_n} \le L,$$

and $|\eta\phi_{n,t}(0)| \le v_n^{-1}$. Hence, for large enough $n$, $\eta\phi_{n,t} \in \Lambda_d(\beta, L)$ for all $\eta \in [-v_n^{-1}, v_n^{-1}]$. By construction, for $s \in [0,1]$ the autoregressive polynomial of $\eta\phi_{n,t}$ is given by $1 - \eta\phi_{n,t}(s)\sum_{i=1}^d u_i z^i$. Since $\lim_{n\to\infty} v_n = \infty$, for any $\rho$, $0 < \rho < 1$, there exists $N$, such that, for all $n \ge N$, $\eta\phi_{n,t} \in \mathcal{S}(\rho)$, $\eta \in [-v_n^{-1}, v_n^{-1}]$, and, thus, $(\eta\phi_{n,t}, \sigma) \in \mathcal{C}$. Using (34) for bounding $\mathbb{E}_{\eta\phi_{n,t},\sigma}[|\mathbf{X}_{k-1,n}|^2]$ in (64), it follows that there exists $M$ depending only on $\rho, \beta, \sigma_+$ and $L$ such that, for all sufficiently large $n$ and for all $\eta \in [-v_n^{-1}, v_n^{-1}]$,

$$\mathcal{I}_n(\eta) \le M\mathcal{I}_\varepsilon \sum_{k\in\mathbb{Z}} \phi^2(kw_n/n - tw_n).$$

Using the fact that $\phi$ is $C^1$ and compactly supported, we have

$$\limsup_{h\to 0}\sup_{x\in\mathbb{R}}h\left|\sum_{k\in\mathbb{Z}}\phi^2(kh-x) - \int \phi^2(t)\,dt\right| = 0.$$

Equation (62) shows that, for large enough $n$ and for all $\eta \in [-v_n^{-1}, v_n^{-1}]$, we have

$$(65) \qquad \mathcal{I}_n(\eta) \le M\mathcal{I}_\varepsilon nw_n^{-1}\int \phi^2(t)\,dt(1 + o(1)).$$

We get that, for all real valued estimators $\hat{\delta}_n := \hat{\delta}_n(\mathbf{X}_{0,n}, X_{1,n}, \ldots, X_{n,n})$, as $n \to \infty$,

$$\sup_{(\boldsymbol{\theta},\sigma)\in\mathcal{C}}\mathbb{E}_{\boldsymbol{\theta},\sigma}[(\hat{\delta}_n - \mathbf{u}^T\boldsymbol{\theta}(t))^2] \ge \int_{-v_n^{-1}}^{v_n^{-1}}v_n\lambda(v_n\eta)\mathbb{E}_{\eta\phi_{n,t},\sigma}[(\hat{\delta}_n - \eta)^2]\,d\eta$$

$$\ge \left(\sup_{\eta\in[-v_n^{-1},v_n^{-1}]}\mathcal{I}_n(\eta) + \mathcal{I}_n(\lambda)\right)^{-1}$$

$$\ge (O(nw_n^{-1} + v_n^2))^{-1},$$

where the first inequality is the Bayesian lower bound of the minimax risk [recall that, for $n$ sufficiently large, $(\eta \phi_{n,t}, \sigma) \in \mathcal{C}$ for all $|\eta| \leq v_n^{-1}$], the second inequality is the so-called van Trees inequality (see [11]) with $\mathcal{I}_n(\eta)$ denoting the Fisher information of the translation model associated with the p.d.f. $v_n \lambda(v_n \cdot)$ and the last inequality is implied by (61) and (65). The proof is concluded by choosing $v_n = w_n^\beta$ and $w_n = n^{1/(1+2\beta)}$.

REMARK 7. Theorem 4 easily extends to cases where the distribution of $\{\mathbf{X}_{0,n}, n \geq 1\}$ depends on $\boldsymbol{\theta}$ in a not too pathological way. Assume, for instance, that $\mathbf{X}_{0,n}$ follows the distribution of a stationary AR process with parameter $(\boldsymbol{\theta}(0), \sigma(0))$ and with a given white noise (see, e.g., [6]). In this case, the lower bound (14) holds for $t > 0$ and $n$ sufficiently large without further assumptions. This clearly follows from the proof: since the distribution of $\mathbf{X}_{0,n}$ depends only on $(\boldsymbol{\theta}(0), \sigma(0))$ and since, for $t > 0$ and $n$ sufficiently large $(\eta \phi_{n,t}(0), \sigma(0)) = (0, \sigma_+)$ does not depend on $\eta$, the computation of $\mathcal{I}_n(\eta)$ applies and the proof proceeds similarly.

**6. Perturbation expansion of the error.** In this section we first derive several approximation results for the error terms $\delta^{(\mathrm{w})}$ and $\delta^{(\mathrm{v})}$ defined by (40) in Section 6.1. Computational estimates are then obtained in Section 6.3 and the proofs of Theorems 6 and 7 are finally given in Section 6.4 and Section 6.5, respectively.

6.1. *General approximation results.* Observe that Theorem 2 only provides a bound of the risk. The approach developed in this section relies upon a perturbation technique (see [1]). Decompose the LHS of (40) as $\delta_{k,n}^{(\bullet)} = J_{k,n}^{(\bullet,0)} + H_{k,n}^{(\bullet,0)}$, with $J_{0,n}^{(\bullet,0)} = 0$, $H_{0,n}^{(\bullet,0)} = 0$ and

$$J_{k+1,n}^{(\bullet,0)} = (I - \mu \mathbb{E}_{\boldsymbol{\theta},\sigma}[\mathrm{F}_\mu(\mathbf{X}_{k,n})]) J_{k,n}^{(\bullet,0)} + \xi_{k,n}^{(\bullet)},$$

$$H_{k+1,n}^{(\bullet,0)} = (I - \mu \mathrm{F}_\mu(\mathbf{X}_{k,n})) H_{k,n}^{(\bullet,0)} + \mu (\mathbb{E}_{\boldsymbol{\theta},\sigma}[\mathrm{F}_\mu(\mathbf{X}_{k,n})] - \mathrm{F}_\mu(\mathbf{X}_{k,n})) J_{k,n}^{(\bullet,0)}.$$

The inhomogeneous first-order difference equation satisfied by $J_{k,n}^{(\bullet,0)}$ yields

$$(66) \qquad J_{k+1,n}^{(\bullet,0)} = \sum_{i=0}^{k} \psi_n(k, i; \mu, \boldsymbol{\theta}, \sigma) \xi_{i,n}^{(\bullet)}, \qquad 0 \leq k < n,$$

where, for all $\mu \geq 0$, $0 \leq i < k \leq n$, and $(\boldsymbol{\theta}, \sigma)$,

$$\psi_n(i, i; \mu, \boldsymbol{\theta}, \sigma) := I \quad \text{and} \quad \psi_n(k, i; \mu, \boldsymbol{\theta}, \sigma) := \prod_{j=i+1}^{k} (I - \mu \mathbb{E}_{\boldsymbol{\theta},\sigma}[\mathrm{F}_\mu(\mathbf{X}_{j,n})]).$$

For ease of notation we write $\mathbf{F}_{k,n} = \mathrm{F}_\mu(\mathbf{X}_{k,n})$, $\overline{\mathbf{F}}_{k,n} := \mathbf{F}_{k,n}(\mu) - \mathbb{E}_{\boldsymbol{\theta},\sigma}[\mathbf{F}_{k,n}(\mu)]$, $\mathbb{E}$ instead of $\mathbb{E}_{\boldsymbol{\theta},\sigma}$, $\Psi_n(k, i)$ instead of $\Psi_n(k, i; \mu)$, $\psi_n(k, i)$ instead of $\psi_n(k, i; \mu,$ $\boldsymbol{\theta}, \sigma$) and so on. This decomposition of the error process $\{\delta_{k,n}^{(\bullet)}, 1 \leq k \leq n\}$ can be extended to further approximation order $s > 0$ as follows:

$$(67) \qquad \delta_{k,n}^{(\bullet)} = J_{k,n}^{(\bullet,0)} + J_{k,n}^{(\bullet,1)} + \cdots + J_{k,n}^{(\bullet,s)} + H_{k,n}^{(\bullet,s)},$$

where

$$J_{k+1,n}^{(\bullet,0)} = (I - \mu \mathbb{E}[\mathbf{F}_{k,n}]) J_{k,n}^{(\bullet,0)} + \xi_{k,n}^{(\bullet)}, \qquad J_{0,n}^{(\bullet,0)} = 0,$$

$$\vdots$$

$$J_{k+1,n}^{(\bullet,r)} = (I - \mu \mathbb{E}[\mathbf{F}_{k,n}]) J_{k,n}^{(\bullet,r)} + \mu \overline{\mathbf{F}}_{k,n} J_{k,n}^{(\bullet,r-1)}, \qquad J_{l,n}^{(\bullet,r)} = 0,$$
$$0 \leq l < r,$$

$$\vdots$$

$$H_{k+1,n}^{(\bullet,s)} = (I - \mu \mathbf{F}_{k,n}) H_{k,n}^{(\bullet,s)} + \mu \overline{\mathbf{F}}_{k,n} J_{k,n}^{(\bullet,s)}, \qquad H_{l,n}^{(\bullet,s)} = 0,$$
$$l = 0, \dots, s.$$

The processes $J_{k,n}^{(\bullet,r)}$ depend linearly on $\xi_{k,n}^{(\bullet)}$ and polynomially in the error $\overline{\mathbf{F}}_{k,n}$. We now show that $J^{(\mathrm{w},0)}$ and $J^{(\mathrm{v},0)}$, respectively, defined by setting $\xi^{(\bullet)} = \xi^{(\mathrm{w})}$ and $\xi^{(\bullet)} = \xi^{(\mathrm{v})}$ in (66), are the main terms in the error terms $\delta^{(\mathrm{w})}$ and $\delta^{(\mathrm{v})}$ defined by (40).

PROPOSITION 20. *Assume* (A1) *with* $q > 4$ *and let* $p \in [1, q/4)$. *Let* $\beta \in (0,1]$, $L > 0$, $0 < \rho < 1$, *and* $0 < \sigma_- \leq \sigma_+$. *Then there exist constants* $M$ *and* $\mu_0 > 0$, *such that, for all* $(\boldsymbol{\theta}, \sigma) \in \mathcal{C}(\beta, L, \rho, \sigma_+, \sigma_-)$, $\mu \in (0, \mu_0]$ *and* $1 \leq k \leq n$,

$$(68) \qquad |J_{k,n}^{(\mathrm{w},0)}| \leq M(n\mu)^{-\beta},$$

$$(69) \qquad \|\delta_{k,n}^{(\mathrm{w})} - J_{k,n}^{(\mathrm{w},0)}\|_{p,\boldsymbol{\theta},\sigma} \leq M\sqrt{\mu}(n\mu)^{-\beta}.$$

PROOF. From (48) in Lemma 19 (which holds under the assumptions of Theorem 15), there exist $\delta > 0$, $\mu_0 > 0$ and $M > 0$ such that, for all $0 \leq i \leq k \leq n$, $\mu \in [0, \mu_0]$ and $(\boldsymbol{\theta}, \sigma) \in \mathcal{C}$,

$$(70) \qquad |\psi_n(k, i; \mu, \boldsymbol{\theta}, \sigma)| \leq M(1 - \delta\mu)^{k-i}.$$

Note that $\psi_n(k, i-1) - \psi_n(k, i) = -\mu \psi_n(k, i) \mathbb{E}[\mathbf{F}_{i,n}]$. As in (54), write, for all $1 \leq k \leq n$,

$$(71) \quad J_{k,n}^{(\mathrm{w},0)} = \psi_n(k-1, 0) \Xi_n^{(\mathrm{w})}(k, 0) + \mu \sum_{j=1}^{k-1} \psi_n(k-1, j) \mathbb{E}[\mathbf{F}_{j-1,n}] \Xi_n^{(\mathrm{w})}(k, j).$$

Using (56), (51) and Lemma C.3 shows (68). By (67) we have $\delta^{(w)} - J^{(w,0)} = H^{(w,0)} = J^{(w,1)} + H^{(w,1)}$, where, for all $1 \le k \le n$,

$$\tag{72} J_{k,n}^{(w,1)} = \mu \sum_{j=0}^{k-1} \psi_n(k-1,j) \overline{\mathbf{F}}_{j,n} J_{j,n}^{(w,0)},$$

$$\tag{73} H_{k,n}^{(w,1)} = \mu \sum_{j=0}^{k-1} \Psi_n(k-1,j) \overline{\mathbf{F}}_{j,n} J_{j,n}^{(w,1)}.$$

Set $\phi_j(\mathbf{x}) = \psi_n(k-1,j) \mathrm{F}_\mu(\mathbf{x}) J_{j,n}^{(w,0)}, j = 0, \ldots, k-1$. Note that, from (49), (70) and (68), for all $\mu \in (0, \mu_0]$, $0 \le j < k \le n$ and $(\boldsymbol{\theta}, \sigma) \in \mathcal{C}$,

$$|\phi_j|_{\mathrm{Li}(1)} \le |\psi_n(k-1,j)| \|J_{j,n}^{(w,0)}\| |\mathrm{F}_\mu|_{\mathrm{Li}(1)} \le M(1-\delta\mu)^{k-j}(\mu n)^{-\beta}.$$

By applying Proposition B.2 componentwise, we get, for all $\mu \in (0, \mu_0]$, $1 \le k \le n$ and $(\boldsymbol{\theta}, \sigma) \in \mathcal{C}$,

$$\left\| \sum_{j=0}^{k-1} \psi_n(k-1,j) \overline{\mathbf{F}}_{j,n} J_{j,n}^{(w,0)} \right\|_{q/2} \le M(\mu n)^{-\beta} \left( \sum_{j=0}^{k-1} (1-\delta\mu)^{k-1-j} \right)^{1/2}.$$

Hence, for all $\mu \in (0, \mu_0]$, $1 \le k \le n$ and $(\boldsymbol{\theta}, \sigma) \in \mathcal{C}$,

$$\tag{74} \|J_{k,n}^{(w,1)}\|_{q/3} \le M\sqrt{\mu}(\mu n)^{-\beta}.$$

Let $u$ be such that $2/q + 2/q + 1/u = 1/p$. Thus, by Theorem 15 and (51), for all $\mu \in (0, \mu_0]$, $1 \le k \le n$ and $(\boldsymbol{\theta}, \sigma) \in \mathcal{C}$,

$$\tag{75} \begin{aligned} \|H_{k,n}^{(w,1)}\|_p &\le \mu \sum_{j=1}^{k-1} \|\Psi_n(k-1,j)\|_u \|\overline{\mathbf{F}}_{j,n}\|_{q/2} \|J_{j,n}^{(w,1)}\|_{q/2} \\ &\le M\mu^{3/2}(\mu n)^{-\beta} \sum_{j=1}^{k-1} (1-\delta\mu)^{k-j} \le M\sqrt{\mu}(\mu n)^{-\beta}. \quad\square \end{aligned}$$

PROPOSITION 21. *Assume* (A1) *with* $q \ge 7$ *and let* $p \in [1, 2q/11)$. *Let* $\beta \in (0,1]$, $L > 0$, $0 < \rho < 1$, *and* $0 < \sigma_- \le \sigma_+$. *Then there exist constants* $M$ *and* $\mu_0$ *such that, for all* $(\boldsymbol{\theta}, \sigma) \in \mathcal{C}(\beta, L, \rho, \sigma_-, \sigma_+)$, $\mu \in (0, \mu_0]$ *and* $1 \le k \le n$,

$$\tag{76} \|J_{k,n}^{(v,0)}(\boldsymbol{\theta}, \sigma)\|_{q, \boldsymbol{\theta}, \sigma} \le M\sqrt{\mu},$$

$$\tag{77} \|\delta_{k,n}^{(v)} - J_{k,n}^{(v,0)}(\boldsymbol{\theta}, \sigma)\|_{p, \boldsymbol{\theta}, \sigma} \le M\mu.$$

PROOF. The Burkholder inequality (see [14], Theorem 2.12) shows that, for all $1 \le k \le n$, $\mu \in (0, \mu_0]$ and $(\boldsymbol{\theta}, \sigma) \in \mathcal{C}$,

$$\|J_{k,n}^{(v,0)}\|_q \le M\mu\sigma_+ \mathrm{L}_q^\star \varepsilon_q^\star \left( \sum_{j=0}^{k-1} |\psi_n(k-1,j)|^2 \right)^{1/2} \le M\mu \left( \sum_{j=0}^{k-1} (1-\delta\mu)^2 \right)^{1/2}$$

and (76) follows from (52) and (70). We now bound

$$J_{k,n}^{(v,1)} = \mu \sum_{j=0}^{k-1} \psi_n(k-1,j)\overline{\mathbf{F}}_{j,n}J_{j,n}^{(v,0)}, \qquad 1 \le k \le n.$$

Let us pick $2 \le k \le n$. By plugging (66) with $\xi^{(\bullet)} = \xi^{(w)}$, we obtain

$$J_{k,n}^{(v,1)} = \mu^2 \sum_{0 \le i < j \le k-1} \phi_j(\mathbf{X}_{j,n})\gamma_{i,j}(\mathbf{X}_{i,n})\sigma_{i+1,n}\varepsilon_{i+1,n}, \tag{78}$$

where, for all $0 \le i < j \le k-1$,

$$\phi_j(\mathbf{x}) := \psi_n(k-1,j)\mathrm{F}_\mu(\mathbf{x}) \quad \text{and} \quad \gamma_{i,j}(\mathbf{x}) := \psi_n(j-1,i)\mathrm{L}_\mu(\mathbf{x}). \tag{79}$$

From (49) and (50) we have, for all $0 \le i < j < k$, $\phi_j \in \mathrm{Li}(\mathbb{R}^d, 1, \mathbb{R}^d \times \mathbb{R}^d; 1)$ and $\gamma_{i,j} \in \mathrm{Li}(\mathbb{R}^d, 1, \mathbb{R}^d; 0)$ and, furthermore, from (70),

$$|\phi_j|_{\mathrm{Li}(1)} \le M(1-\delta\mu)^{k-j} \quad \text{and} \quad |\gamma_{i,j}|_{\mathrm{Li}(0)} \le M(1-\delta\mu)^{j-i},$$

where, as usual, $M$ and $\delta$ are positive constants, not depending on indices $i, j, k, n$, on $\mu \in [0, \mu_0]$ or on $(\boldsymbol{\theta}, \sigma) \in \mathcal{C}$. The following uniform bounds follow:

$$\sup_{0 < j < t} |\phi_j|_{\mathrm{Li}(1)} \le M, \qquad\qquad \sum_{j=1}^{t-1} |\phi_j|_{\mathrm{Li}(1)} \le M\mu^{-1},$$

$$\sup_{0 \le i < j < t} |\gamma_{i,j}|_{\mathrm{Li}(0)} \le M, \qquad \sup_{0 < j < t}\left(\sum_{0 \le i < j} |\gamma_{i,j}|_{\mathrm{Li}(0)}^2\right)^{1/2} \le M\mu^{-1/2}.$$

By applying Proposition B.3 componentwise, we obtain $\|J_{k,n}^{(v,1)}\|_{2q/7} \le M\mu$ uniformly over $1 \le k \le n$, $\mu \in [0, \mu_0]$ and $(\boldsymbol{\theta}, \sigma) \in \mathcal{C}$. As in (75), let $u > 1$ be such that $u^{-1} + 2/q + 7/2q = 1/p$. Then, for all $1 \le k \le n$, $\mu \in [0, \mu_0]$ and $(\boldsymbol{\theta}, \sigma) \in \mathcal{C}$,

$$\|H_{k,n}^{(v,1)}\|_p \le \mu \sum_{j=1}^{k-1} \|\Psi_n(k-1,j)\|_u \|\overline{\mathbf{F}}_{j,n}\|_{q/2}\|J_{j,n}^{(v,1)}\|_{2q/7}$$

and thus, $\|H_{k,n}^{(v,1)}\|_p \le M\mu^2\mathbf{F}_{q/2}^\star\sum_{j=1}^{k-1}(1-\delta\mu)^{k-j}$, which yields (77). $\square$

6.2. *Proof of Proposition 5.* We may write, for all $1 \le k \le n$,

$$\mathbf{X}_{k,n} = \beta_n(k, 0; \boldsymbol{\theta})\mathbf{X}_{0,n} + \sum_{j=1}^{k} \beta_n(k, j; \boldsymbol{\theta})\boldsymbol{\sigma}_{j,n}\varepsilon_{j,n},$$

where $\boldsymbol{\sigma}$ and $\beta_n$ are defined right after (30) and in (33), respectively. Thus,

$$\mathbb{E}_{\boldsymbol{\theta},\sigma}[\mathbf{X}_{k,n}\mathbf{X}_{k,n}^T] = \beta_n(k,0)\mathbb{E}[\mathbf{X}_{0,n}\mathbf{X}_{0,n}^T]\beta_n(k,0)^T$$

$$+ \sum_{l=0}^{k-1} \beta_n(k, k-l)\boldsymbol{\sigma}_{k-l,n}(\beta_n(k, k-l)\boldsymbol{\sigma}_{k-l,n})^T.$$

Let $\boldsymbol{\theta} \in \mathcal{S}(\rho)$ with $\rho > 1$, $t \in [0,1]$ and let $\{Z_k, k \in \mathbb{Z}\}$ denote the stationary AR($\boldsymbol{\theta}(t), \sigma(t)$) with i.i.d. centered unit variance innovations denoted by $\{\varepsilon_k\}_{k \in \mathbb{Z}}$. Recall that $\Sigma(t, \boldsymbol{\theta}, \sigma)$ denotes the $d \times d$ covariance matrix of $\{Z_k, k \in \mathbb{Z}\}$. Then, using classical results on AR models (see [4]), we have $[Z_k \cdots Z_{k-d}]^T = \sum_{l \geq 0} \Theta^l(t, \boldsymbol{\theta}) \boldsymbol{\sigma}(t) \varepsilon_{k-l}$, where $\boldsymbol{\sigma}(t) := [\sigma(t)\, 0 \cdots 0]^T$, $\Theta$ is defined by (29) and the convergence holds in the $L^2$ sense. It follows that

$$\Sigma(t, \boldsymbol{\theta}, \sigma) = \sum_{l=0}^{\infty} \Theta^l(t, \boldsymbol{\theta}) \boldsymbol{\sigma}(t) (\Theta^l(t, \boldsymbol{\theta}) \boldsymbol{\sigma}(t))^T.$$

Denote $\Theta_{k,n} := \Theta(k/n, \boldsymbol{\theta})$ and $\Sigma_{k,n} := \Sigma(k/n, \boldsymbol{\theta}, \sigma)$. We obtain

$$\begin{aligned}
&\mathbb{E}_{\boldsymbol{\theta}, \sigma}[\mathbf{X}_{k,n} \mathbf{X}_{k,n}^T] - \Sigma_{k,n} \\
(80) \quad &= \beta_n(k,0)(\mathbb{E}[\mathbf{X}_{0,n} \mathbf{X}_{0,n}^T] - \Sigma_{0,n}) \beta_n(k,0)^T \\
&\quad + \sum_{l=0}^{k-1} (\beta_n(k,k-l) \boldsymbol{\sigma}_{k-l,n} (\beta_n(k,k-l) \boldsymbol{\sigma}_{k-l,n})^T - \Theta_{k,n}^l \boldsymbol{\sigma}_{k,n} (\Theta_{k,n}^l \boldsymbol{\sigma}_{k,n})^T) \\
&\quad + \sum_{l=k}^{\infty} (\beta_n(k,0) \Theta_{0,n}^{l-k} \boldsymbol{\sigma}_{0,n} (\beta_n(k,0) \Theta_{0,n}^{l-k} \boldsymbol{\sigma}_{0,n})^T - \Theta_{k,n}^l \boldsymbol{\sigma}_{k,n} (\Theta_{k,n}^l \boldsymbol{\sigma}_{k,n})^T).
\end{aligned}$$

Note that, for any matrices $A_1, \ldots, A_r$ and $B_1, \ldots, B_r$ with compatible sizes,

$$(81) \quad \prod_{i=1}^{r} A_i - \prod_{i=1}^{r} B_i = \sum_{j=1}^{r} \left( \prod_{k=1}^{j-1} A_k \right) (A_j - B_j) \left( \prod_{k=j+1}^{r} B_k \right).$$

From uniform exponential stability (see Proposition 13 and its proof), there exists $M$ such that, for all $1 \leq l \leq k \leq n$ and $(\boldsymbol{\theta}, \sigma) \in \mathcal{C}^\star$, $|\beta_n(k,k-l)| \leq M\tau^l$ and $|\Theta_{k,n}^l| \leq M\tau^l$ and, thus,

$$|\beta_n(k,k-l) - \Theta_{k,n}^l| \leq M\tau^l \sum_{j=0}^{l-1} |\Theta_{k-j,n} - \Theta_{k,n}| \leq M n^{-\beta} \tau^l l^{\beta+1}.$$

Similarly, there exists $M$ such that, for all $1 \leq l \leq k \leq n$ and $(\boldsymbol{\theta}, \sigma) \in \mathcal{C}^\star$,

$$\begin{aligned}
&|\Theta_{0,n}^{l-k} - \Theta_{k,n}^{l-k}| \\
&\quad \leq M n^{-\beta} \tau^{l-k} (l-k) k^\beta, \\
&|\beta_n(k,k-l) \boldsymbol{\sigma}_{k-l,n} (\beta_n(k,k-l) \boldsymbol{\sigma}_{k-l,n})^T - \Theta_{k,n}^l \boldsymbol{\sigma}_{k,n} (\Theta_{k,n}^l \boldsymbol{\sigma}_{k,n})^T| \\
&\quad \leq M n^{-\beta} \tau^{2l} l^{\beta+1}, \\
&|\beta_n(k,0) \Theta_{0,n}^{l-k} \boldsymbol{\sigma}_{0,n} (\beta_n(k,0) \Theta_{0,n}^{l-k} \boldsymbol{\sigma}_{0,n})^T - \Theta_{k,n}^l \boldsymbol{\sigma}_{k,n} (\Theta_{k,n}^l \boldsymbol{\sigma}_{k,n})^T| \\
&\quad \leq M n^{-\beta} \tau^{2l} l k^\beta.
\end{aligned}$$

The result follows by inserting these bounds in (80). $\quad \square$

6.3. *Further approximation results.* To derive tractable asymptotic risk estimates, we need to derive approximate expressions for $J_{k,n}^{(w,0)}$ and $J_{k,n}^{(v,0)}$. We first derive approximations for $\mathbb{E}[\mathbf{F}_{k,n}], 1 \leq k \leq n$, and related quantities.

LEMMA 22. *Assume* (A1) *with* $q \geq 4$. *Let* $\beta \in (0,1]$, $L > 0$, $0 < \rho < \tau < 1$, *and* $0 < \sigma_- \leq \sigma_+$. *Then there exist positive constants* $\delta, \nu_0$ *and* $M$ *such that, for all* $0 \leq k \leq l \leq n$, $\nu \in [0, \nu_0]$ *and* $(\boldsymbol{\theta}, \sigma) \in \mathcal{C}^\star(\beta, L, \rho, \sigma_-, \sigma_+)$,

$$\tag{82} |(I - \nu\Sigma(l/n, \boldsymbol{\theta}, \sigma))^{l-k}| \leq (1 - \delta\nu)^{l-k},$$

$$\tag{83} |\Sigma(l/n, \boldsymbol{\theta}, \sigma) - \mathbb{E}_{\boldsymbol{\theta}, \sigma}[\mathbf{F}_\nu(\mathbf{X}_{k,n})]| \leq M(\tau^k + n^{-\beta}(l-k+1)^\beta + \nu),$$

$$\tag{84} |\Sigma(l/n, \boldsymbol{\theta}, \sigma) - \mathbb{E}_{\boldsymbol{\theta}, \sigma}[\mathbf{L}_\nu\mathbf{L}_\nu^T(\mathbf{X}_{k,n})]| \leq M(\tau^k + n^{-\beta}(l-k+1)^\beta + \nu).$$

PROOF. By continuity of $(\boldsymbol{\theta}, z, t) \mapsto |\theta(z;t)|$ [see (6)] and since

$$\left\{ \boldsymbol{\vartheta} \in \mathbb{R}^d : |\boldsymbol{\vartheta}| \leq L, 1 - \sum_{i=1}^d \boldsymbol{\vartheta}_i z^i \neq 0 \text{ for all } |z| \leq \rho^{-1} \right\}$$

is a compact set, there exist $\delta > 0$ and $M > 0$, such that for all $(\boldsymbol{\theta}, \sigma) \in \mathcal{C}^\star$,

$$\tag{85} \begin{aligned} \delta &\leq \inf_{|z|=1} \inf_{t \in [0,1]} \frac{\sigma^2(t)}{|\theta(z;t)|^2} \leq \lambda_{\min}(\Sigma(t; \boldsymbol{\theta}, \sigma)) \\ &\leq \lambda_{\max}(\Sigma(t; \boldsymbol{\theta}, \sigma)) \leq \sup_{|z|=1} \sup_{t \in [0,1]} \frac{\sigma^2(t)}{|\theta(z;t)|^2} \leq M. \end{aligned}$$

Equation (82) then follows from Lemma C.1. Similarly, there exists $M < \infty$, such that for all $(\boldsymbol{\theta}, \sigma) \in \mathcal{C}^\star$ and all $0 \leq s \leq t \leq 1$,

$$\tag{86} |\Sigma(t; \boldsymbol{\theta}, \sigma) - \Sigma(s; \boldsymbol{\theta}, \sigma)| \leq M(t-s)^\beta.$$

By Proposition 5 we get, for all $1 \leq k \leq l \leq n$ and $(\boldsymbol{\theta}, \sigma) \in \mathcal{C}^\star$,

$$\tag{87} \begin{aligned} &|\Sigma(l/n; \boldsymbol{\theta}, \sigma) - \mathbb{E}[\mathbf{X}_{k,n}\mathbf{X}_{k,n}^T]| \\ &\quad \leq |\Sigma(l/n; \boldsymbol{\theta}, \sigma) - \Sigma(k/n; \boldsymbol{\theta}, \sigma)| + |\Sigma(k/n; \boldsymbol{\theta}, \sigma) - \mathbb{E}[\mathbf{X}_{k,n}\mathbf{X}_{k,n}^T]| \\ &\quad \leq M(\tau^k + n^{-\beta}(l-k+1)^\beta). \end{aligned}$$

This is (83) and (84) with $\nu = 0$. One easily shows that, for all $\mathbf{x} \in \mathbb{R}^d$ and $\nu \geq 0$, $|\mathbf{F}_\nu(\mathbf{x}) - \mathbf{x}\mathbf{x}^T| \leq \nu|\mathbf{x}|^4$ and $|\mathbf{L}_\nu\mathbf{L}_\nu^T(\mathbf{x})| \leq 2\nu|\mathbf{x}|^4$. Since $q \geq 4$, we deduce (83) and (84) for $\nu > 0$ from (87) and uniform $L^4$ boundedness. $\square$

Let $\rho \in (0,1)$ and let $\boldsymbol{\theta} \in \mathcal{S}(\rho)$ and $\sigma : [0,1] \to \mathbb{R}^+$. Define the following sequence of recurrence equations applying to some increment process $\{\xi_{k,n}^{(\bullet)}, 0 \leq k \leq n\}$:

$$\tag{88} \tilde{J}_{k+1,n}^{(\bullet)}(\boldsymbol{\theta}, \sigma) := \sum_{j=0}^k (I - \mu\Sigma((k+1)/n, \boldsymbol{\theta}, \sigma))^{k-j}\xi_{j,n}^{(\bullet)}, \qquad 0 \leq k \leq n.$$

We now show that $J^{(w,0)}$ and $J^{(v,0)}$ may be approximated by $\tilde{J}^{(w)}$ and $\tilde{J}^{(v)}$, respectively, defined by setting $\xi^{(\bullet)} = \xi^{(w)}$ and $\xi^{(\bullet)} = \xi^{(v)}$ in (88) and then we compute asymptotic equivalents of $\tilde{J}^{(w)}$ and of $\tilde{J}^{(v)}$'s variance, respectively.

PROPOSITION 23. *Assume* (A1) *with* $q \geq 4$. *Let* $\beta \in (0,1]$, $L > 0$, $\rho < 1$, *and* $0 < \sigma_- \leq \sigma_+$. *Then there exist constants* $M > 0$ *and* $\mu_0 > 0$ *such that, for all* $(\boldsymbol{\theta}, \sigma) \in \mathcal{C}^\star(\beta, L, \rho, \sigma_-, \sigma_+)$, $\mu \in (0, \mu_0]$ *and* $1 \leq k \leq n$,

$$(89) \qquad |J_{k,n}^{(w,0)}(\boldsymbol{\theta}, \sigma) - \tilde{J}_{k,n}^{(w)}(\boldsymbol{\theta}, \sigma)| \leq M(\mu n)^{-\beta}((\mu n)^{-\beta} + \mu),$$

$$(90) \qquad \|J_{k,n}^{(v,0)}(\boldsymbol{\theta}, \sigma) - \tilde{J}_{k,n}^{(v)}(\boldsymbol{\theta}, \sigma)\|_{q,\boldsymbol{\theta},\sigma} \leq M\sqrt{\mu}((\mu n)^{-\beta} + \mu).$$

PROOF. Let $\Delta_n(k, j; \mu, \boldsymbol{\theta}, \sigma) := \psi_n(k, j; \mu, \boldsymbol{\theta}, \sigma) - (I - \mu\Sigma(k/n, \boldsymbol{\theta}, \sigma))^{k-1-j}$ for all $0 \leq j \leq k \leq n$, $\mu \geq 0$ and $(\boldsymbol{\theta}, \sigma)$. The definitions of $J^{(\bullet,0)}$ and $\tilde{J}^{(\bullet)}$ yield

$$(91) \qquad J_{k,n}^{(\bullet,0)} - \tilde{J}_{k,n}^{(\bullet)} = \sum_{j=0}^{k-1} \Delta_n(k-1, j)\xi_{j,n}^{(\bullet)}, \qquad 1 \leq k \leq n.$$

From (81) we get, for all $0 \leq j < k \leq n$, $\mu \geq 0$ and $(\boldsymbol{\theta}, \sigma) \in \mathcal{C}^\star$,

$$\Delta_n(k-1, j) = \mu \sum_{i=j}^{k-1} \psi_n(k-1, i+1)(\mathbb{E}[\mathbf{F}_{i,n}] - \Sigma(k/n))(I - \mu\Sigma(k/n))^{i-j-1}.$$

Using (70), (82) and (83), there exist $\delta > 0$, $\mu_0 > 0$ and $M > 0$, such that, for all $1 \leq j < k \leq n$, $\mu \in (0, \mu_0]$ and $(\boldsymbol{\theta}, \sigma) \in \mathcal{C}^\star$,

$$(92) \quad |\Delta_n(k-1, j)| \leq M\mu(1 - \delta\mu)^{k-1-j}(\tau^j + n^{-\beta}(k-j)^{\beta+1} + \mu(k-j)).$$

We further write, for all $1 \leq j < k \leq n$,

$$\begin{aligned}
\Delta_n(k-1&, j) - \Delta_n(k-1, j-1) \\
&= \mu\Delta_n(k-1, j)\mathbb{E}[\mathbf{F}_{j-1,n}] \\
&\quad + \mu(I - \mu\Sigma(k/n))^{k-j-1}(\mathbb{E}[\mathbf{F}_{j-1,n}] - \Sigma(k/n)).
\end{aligned}$$

Applying (92), (82) and (83) and observing that $\mathbf{F}_1^\star$ is finite, we get that there exist $\delta > 0$, $\mu_0 > 0$, $\tau \in (\rho, 1)$ and $M > 0$, such that, for all $1 \leq j < k \leq n$, $\mu \in (0, \mu_0]$ and $(\boldsymbol{\theta}, \sigma) \in \mathcal{C}^\star$,

$$\begin{aligned}
(93) \qquad & |\Delta_n(k-1, j) - \Delta_n(k-1, j-1)| \\
& \leq M\mu(1-\delta\mu)^{k-j-1}(n^{-\beta}(k-j)^\beta(\mu(k-j)+1) \\
& \quad + \mu(\mu(k-j)+1) + \tau^j).
\end{aligned}$$

By integrating (91) by parts, for all $1 \leq k \leq n$, $J_{k,n}^{(\bullet,0)} - \tilde{J}_{k,n}^{(\bullet)}$ reads

$$\Delta_n(k-1,0)\Xi_n^{(\bullet)}(k,0) + \sum_{j=1}^{k-1}(\Delta_n(k-1,j) - \Delta_n(k-1,j-1))\Xi_n^{(\bullet)}(k,j).$$

Using (56) and (58) to bound $|\Xi_n^{(w)}|$ and $\|\Xi_n^{(v)}\|_q$ respectively, together with (92), (93) and Lemma C.3, there exist $\delta > 0$, $\mu_0 > 0$ and $M$ such that, for all $1 \leq k \leq n$, $\mu \in (0, \mu_0]$ and $(\boldsymbol{\theta}, \sigma) \in \mathcal{C}$,

$$|J_{k,n}^{(w,0)} - \tilde{J}_{k,n}^{(w)}| \leq M\left(\frac{(\mu n)^{-\beta} + \mu}{(\mu n)^\beta} + \mu\sum_{j=0}^{k-1}(1-\delta\mu)^j\tau^{k-j}\left(\frac{j+1}{n}\right)^\beta\right),$$

$$\|J_{k,n}^{(v,0)} - \tilde{J}_{k,n}^{(v)}\|_q \leq M\left(\sqrt{\mu}((\mu n)^{-\beta} + \mu) + \mu^2\sum_{j=0}^{k-1}(1-\delta\mu)^j\tau^{k-j}\sqrt{j+1}\right).$$

Using Lemma C.3, we have, for any $\alpha \geq 0$,

$$\sum_{j=0}^{k-1}(1-\delta\mu)^j\tau^{k-j}(j+1)^\alpha$$

$$\leq (1-\tau)^{-1}\sup_{j \in \mathbb{N}}(1-\delta\mu)^j(j+1)^\alpha$$

$$\leq C(1-\tau)^{-1}(\delta\mu)^{-\alpha}.$$

We thus obtain (89) and (90). $\quad\square$

PROPOSITION 24. *Let $\beta \in (0,1]$, $L > 0$, $\rho < 1$, and $0 < \sigma_- \leq \sigma_+$ and let $(\boldsymbol{\theta}, \sigma) \in \mathcal{C}^\star(\beta, L, \rho, \sigma_-, \sigma_+)$. Let $t \in (0,1]$ and assume that there exists $\boldsymbol{\theta}_{t,\beta} \in \mathbb{R}^d$, $L' > 0$ and $\beta' > \beta$ such that (18) holds for all $u \in [0,t]$. Then for all $\mu \in (0, \mu_0]$ and $n \geq 1$,*

$$\left|\tilde{J}_{[tn],n}^{(w)} - \frac{\Gamma(\beta+1)}{(\mu n)^\beta}\Sigma^{-\beta}(t)\boldsymbol{\theta}_{t,\beta}\right| \leq M\left((\mu n)^{-\beta'} + \frac{n^{-\beta} + (1-\delta\mu)^{tn}}{(\mu n)^\beta}\right),$$

*where $M$ and $\mu_0$ are positive constants depending only on $\beta, L, \rho, \sigma_-, \sigma_+, L'$ and $\beta'$.*

PROOF. Let us write

$$\tilde{J}_{[tn],n}^{(w)} - J_n$$

$$= \sum_{j=0}^{[tn]-1}(I - \mu\Sigma([tn]/n))^{[tn]-1-j}$$

$$\times [\xi_{j,n}^{(w)} - n^{-\beta}(([tn]-j)^\beta - ([tn]-1-j)^\beta)\boldsymbol{\theta}_{t,\beta}]$$

where, within this proof section, we denote

$$J_n := n^{-\beta} \sum_{j=0}^{[tn]-1} (I - \mu\Sigma([tn]/n))^{[tn]-1-j}(([tn]-j)^\beta - ([tn]-1-j)^\beta)\boldsymbol{\theta}_{t,\beta}.$$

Using (18), the partial sums of the terms within brackets in the next-to-last equation satisfy, for all $0 \leq j \leq [tn] - 1$,

$$\left| \sum_{i=j}^{[tn]-1} [\xi_{i,n}^{(w)} - n^{-\beta}(([tn]-i)^\beta - ([tn]-1-i)^\beta)\boldsymbol{\theta}_{t,\beta}] \right|$$

$$= |\boldsymbol{\theta}(j/n) - \boldsymbol{\theta}([tn]/n) + n^{-\beta}[tn]-j)^\beta\boldsymbol{\theta}_{t,\beta}|$$

$$\leq L'n^{-\beta'}([tn]-j)^{\beta'}.$$

Integration by parts with this bound and (82), and then Lemma C.3, give that there exists a constant $M$ such that, for all $\mu \in (0, \mu_0]$ and $n \geq 1$,

$$|\tilde{J}_{[tn],n}^{(w)} - J_n| \leq M(\mu n)^{-\beta'}.$$

Now, from (82) and using Lemma C.3, we have, for all $\mu \in (0, \mu_0]$ and $n \geq 1$,

$$|J_n - n^{-\beta}S_\beta(I - \mu\Sigma([tn]/n))\boldsymbol{\theta}_{t,\beta}|$$

$$\leq Mn^{-\beta}(1-\delta\mu)^{[tn]} \sum_{l\geq 1}(1-\delta\mu)^l l^{\beta-1}$$

$$\leq M(1-\delta\mu)^{[tn]}(\mu n)^{-\beta},$$

where $S_\beta(A) := \sum_{i=0}^\infty A^i((i+1)^\beta - i^\beta)$. Using the fact that $A \to S_\beta(A)$ is a power series with unit radius of convergence, (82) and (86), from the mean value theorem and Lemma C.3, we have, for all $\mu \in (0, \mu_0]$ and all $n \geq 1$,

$$|S_\beta(I - \mu\Sigma([tn]/n)) - S_\beta(I - \mu\Sigma(t))| \leq M\mu n^{-\beta} \sum_{i\geq 1}(1-\delta\mu)^i i^\beta$$

$$\leq M(\mu n)^{-\beta}.$$

Collecting the last three inequalities together with Lemma C.3, we obtain the result. □

PROPOSITION 25. *Assume* (A1) *with* $q \geq 4$. *Let* $\beta \in (0,1]$, $L > 0$, $\rho < 1$, *and* $0 < \sigma_- \leq \sigma_+$. *Then there exist constants* $M > 0$ *and* $\mu_0 > 0$ *such that, for all* $(\boldsymbol{\theta}, \sigma) \in \mathcal{C}^\star(\beta, L, \rho, \sigma_-, \sigma_+)$, $\mu \in (0, \mu_0]$ *and* $1 \leq k \leq n$,

$$\left| \mathbb{E}_{\boldsymbol{\theta},\sigma}[\tilde{J}_{k,n}^{(v)}\tilde{J}_{k,n}^{(v)T}] - \mu\frac{\sigma^2(k/n)}{2}I \right| \leq M\mu(\mu + (\mu n)^{-\beta} + (1-\delta\mu)^k).$$

PROOF. Since $\{\xi^{(v)}{}_{i,n}, j \geq 0\}$ is a martingale increment sequence, for any $1 \leq k \leq n$ $\mathbb{E}[\tilde{J}^{(v)}_{k,n} \tilde{J}^{(v)T}_{k,n}]$ reads

$$\mu^2 \sum_{j=0}^{k-1} (I - \mu\Sigma(k/n))^{k-1-j} \sigma^2_{j,n} \mathbb{E}[L_\mu(\mathbf{X}_{j,n}) L_\mu^T(\mathbf{X}_{j,n})] (I - \mu\Sigma(k/n))^{k-1-j}$$

$$= \mu\sigma^2_{k,n}(G_{k,n} - \tilde{G}_{k,n}) + R_{k,n},$$

where, for all $1 \leq k \leq n$,

$$G_{k,n} := \mu \sum_{j=-\infty}^{k-1} (I - \mu\Sigma(k/n))^{k-1-j} \Sigma(k/n)(I - \mu\Sigma(k/n))^{k-1-j},$$

$$\tilde{G}_{k,n} := \mu \sum_{j=-\infty}^{-1} (I - \mu\Sigma(k/n))^{k-1-j} \Sigma(k/n)(I - \mu\Sigma(k/n))^{k-1-j},$$

$$|R_{k,n}| \leq M\mu^2 \sum_{j=0}^{k-1} (1-\delta\mu)^{2(k-1-j)}(\tau^j + n^{-\beta}(k-j)^\beta + \mu)$$

$$\leq M\mu(\mu + (\mu n)^{-\beta}).$$

For bounding $R_{k,n}$, we have used (82), (84), $\sigma \in \Lambda(\beta, L)$ and then Lemma C.3. From (82) and (85), we have, for all $\mu \in (0, \mu_0]$, $1 \leq k \leq n$,

(94) $$|\tilde{G}_{k,n}| \leq M(1-\delta\mu)^k.$$

From the previous bounds, we obtain, for all $\mu \in (0, \mu_0]$ and $1 \leq k \leq n$,

(95) $$|\mathbb{E}[\tilde{J}^{(v)}_{k,n} \tilde{J}^{(v)T}_{k,n}] - \mu\sigma^2_{k,n} G_{k,n}| \leq M\mu(\mu + (\mu n)^{-\beta} + (1-\delta\mu)^k).$$

Now, by definition of $G_{k,n}$ we have

$$(I - \mu\Sigma(k/n))G_{k,n}(I - \mu\Sigma(k/n)) + \mu\Sigma(k/n) = G_{k,n}, \qquad 1 \leq k \leq n,$$

which gives, for all $1 \leq k \leq n$,

$$(\Sigma(k/n)(I - 2G_{k,n}) + (I - 2G_{k,n})\Sigma(k/n)) = 2\mu\Sigma(k/n)G_{k,n}\Sigma(k/n).$$

From (94) $\sup_{1 \leq k \leq n} |G_{k,n}| < \infty$ uniformly over $\mu \in (0, \mu_0]$. We thus have, for all $\mu \in (0, \mu_0]$ and $1 \leq k \leq n$,

$$|\Sigma(k/n)(I - 2G_{k,n}) + (I - 2G_{k,n})\Sigma(k/n)| \leq M\mu.$$

For any $d \times d$ matrix $C$ and positive definite matrix $S$, the equation $SB + BS = C$ has a unique solution $B$ linear in $C$ and continuous in $S$ over the set of positive definite matrices (see [16], Corollary 4.4.10). Hence, we obtain, for all $\mu \in (0, \mu_0]$ and $1 \leq k \leq n$, $|I - 2G_{k,n}(\mu)| \leq M\mu$, which with (95) gives the claimed bound. $\square$

6.4. *Proof of Theorem* 6. We use the decomposition of $\delta$ as

$$\delta^{(\mathrm{u})} + (\delta^{(\mathrm{w})} - J^{(\mathrm{w},0)}) + (\delta^{(\mathrm{v})} - J^{(\mathrm{v},0)})$$

$$+ (J^{(\mathrm{w},0)} - \tilde{J}^{(\mathrm{w})}) + (J^{(\mathrm{v},0)} - \tilde{J}^{(\mathrm{v})}) + \tilde{J}^{(\mathrm{w})} + \tilde{J}^{(\mathrm{v})}.$$

Let $\eta \in (0,1)$. Applying (53), (69), (77), (89) and (90), there exists $M > 0$ such that, for all $t \in [\eta, 1]$, $\mu \in (0, \mu_0]$ and $n \geq 1$,

$$\|\hat{\boldsymbol{\theta}}_n(t; \mu) - \boldsymbol{\theta}(t) - \tilde{J}^{(\mathrm{w})}_{[tn], n} - \tilde{J}^{(\mathrm{v})}_{[tn], n}\|_{2, \boldsymbol{\theta}, \sigma}$$

$$\leq M(\sqrt{\mu}(\mu n)^{-\beta} + (\mu n)^{-2\beta} + \mu).$$

We then obtain (19) by applying Proposition 24 and Proposition 25, and using the fact that $\sigma$ is $\beta$-Lipschitz to approximate $\sigma^2([tn]/n)$ by $\sigma^2(t)$.

6.5. *Proof of Theorem* 7. We use $\delta = \delta^{(\mathrm{u})} + \delta^{(\mathrm{v})} + (\delta^{(\mathrm{w})} - J^{(\mathrm{w},0)}) + (J^{(\mathrm{w},0)} - \tilde{J}^{(\mathrm{w})}) + \tilde{J}^{(\mathrm{w})}$. Observe that there exists $C > 0$ such that $\Lambda(\beta, L) \subseteq \Lambda(1, CL)$. Hence, we may apply (53), (59), (69) and (89), so that there exist $M$ and $\mu_0$ such that, for all $\mu \in (0, \mu_0]$,

$$(96) \qquad \sup_{(\boldsymbol{\theta}, \sigma) \in \mathcal{C}^\star} \sup_{t \in [\eta, 1]} \|\hat{\boldsymbol{\theta}}_n(t; \mu) - \boldsymbol{\theta}(t) - \tilde{J}^{(\mathrm{w})}_{[tn], n}\|_{p, \boldsymbol{\theta}, \sigma} \leq M(\sqrt{\mu} + (\mu n)^{-2}).$$

Now, since $\beta > 1$, for all $(\boldsymbol{\theta}, \sigma) \in \mathcal{C}^\star$ and $t \in (0,1]$, we may apply the Taylor expansion $\boldsymbol{\theta}(u) = \boldsymbol{\theta}(t) + \dot{\boldsymbol{\theta}}(v)(u - t)$, where $v \in [u, t]$, which yields $|\boldsymbol{\theta}(u) - \boldsymbol{\theta}(t) + \dot{\boldsymbol{\theta}}(t)(t - u)| \leq |\dot{\boldsymbol{\theta}}(v) - \dot{\boldsymbol{\theta}}(t)||u - t| \leq L|t - u|^{\beta - 1}$. Hence (18) holds ($\beta = 1$ and $\beta'$ equal to the actual $\beta$) at every point $t > 0$ and we may apply Proposition 24 for computing $\tilde{J}^{(\mathrm{w})}$, which easily yields the result.

## APPENDIX A

PROOF OF THEOREM 16. We first derive two simple lemmas valid under the assumptions of Theorem 16. We let $\nu_0$ and $\delta$ denote some constants depending only on $r$, $R_1$, $\nu_1$, $\alpha_1$, $C_1$, $\lambda$, $B$ and $s$ and we write $\mathbb{E}_k$ for $\mathbb{E}^{\mathcal{F}_k}$.

LEMMA A.1. *For any $a \geq 1$, there exist $\delta > 0$ and $\nu_0 > 0$ such that, for all $k \in \mathbb{N}$ and $\nu \in [0, \nu_0]$,*

$$(97) \qquad \mathbf{I}(\phi_k \leq R_1) \mathbb{E}_k \left\{ \left| \prod_{i=k+1}^{k+r} (I - \nu A_i) \right|^a \right\} \leq e^{-\delta \nu}.$$

PROOF. Under (C-1) we have $|I - \nu A_k| \leq 1$ for all $k \in \mathbb{N}$ and $\nu \in [0, \nu_1]$ (see Lemma C.1) so that we may assume $a = 1$ without loss of generality. We write

$$(98) \qquad \prod_{l=k+1}^{k+r} (I - \nu A_l) = I - \nu D_k + S_k,$$

where

$$D_k := \sum_{l=k+1}^{k+r} A_l \quad \text{and} \quad S_k := \sum_{j=2}^{r} (-1)^j \nu^j \sum_{1 \le i_1 < \cdots < i_j \le r} A_{k+i_j} \ldots A_{k+i_1}.$$

For $\beta \in (1/s, 1)$, where $s$ is defined in (C-3) and $\nu \ge 0$, denote

$$\mathsf{B}_k(\nu) := \{|A_{k+1}| \le \nu^{-\beta}, \ldots, |A_{k+r}| \le \nu^{-\beta}\}$$

and $\mathsf{B}_k^c(\nu)$ its complementary set. From (C-1), we have that, for all $\nu \in (0, \nu_1]$, $|\nu D_k| \le r\nu^{1-\beta} \mathbf{I}(\mathsf{B}_k(\nu)) + r\mathbf{I}(\mathsf{B}_k^c(\nu))$. Choosing $\nu_2 \in (0, \nu_1]$ such that $r\nu_2^{1-\beta} \le 1$, we get that, for all $\nu \in [0, \nu_2]$, $2|\nu D_k|\mathbf{I}(|\nu D_k| > 1) \le 2r\mathbf{I}(\mathsf{B}_k^c(\nu))$. Hence, using (98) and Lemma C.1, we obtain, for all $\nu \in [0, \nu_2]$,

$$\left| \prod_{l=k+1}^{k+r} (I - \nu A_l) \right| \le 1 - \nu \lambda_{\min}(D_k) + 2r\mathbf{I}(\mathsf{B}_k^c(\nu)) + |S_k|.$$

Equation (97) easily follows from this bound with (C-2) and the two following inequalities, which will be shown to hold for all $\nu \in [0, \nu_1]$:

$$(99) \qquad \mathbf{I}(\phi_k \le R_1)\mathbb{E}_k[\mathbf{I}(\mathsf{B}_k^c(\nu))] \le C_1\nu^{s\beta},$$

$$(100) \qquad \mathbf{I}(\phi_k \le R_1)\mathbb{E}_k[|S_k|] \le MC_1^{(s \wedge 2)/s} \nu^{s \wedge 2},$$

where $M$ is some constant depending only on $r$. We now conclude the proof by showing these two last inequalities successively.

Using the Markov inequality, we obtain $\mathbb{E}_k[\mathbf{I}(\mathsf{B}_k^c(\nu))] \le \sum_{l=k+1}^{k+r} \mathbb{P}_k\{|A_l| > \nu^{-\beta}\} \le \nu^{s\beta} \sum_{l=k+1}^{k+r} \mathbb{E}_k[|A_l|^s]$, which implies (99) using (C-3).

For all $j = 2, \ldots, r$, for all ordered $j$-tuples $1 \le i_1 < \cdots < i_j \le r$, using (C-1), we have, for all $\nu \in [0, \nu_1]$, $\nu^j|A_{k+i_1} \ldots A_{k+i_j}| \le \nu^2|A_{k+i_1} A_{k+i_2}|$. Hence, for some constant $M_1$ depending only on $r$, for all $\nu \in [0, \nu_1]$,

$$\mathbb{E}_k[|S_k|] \le M_1\nu^2 \sup_{1 \le i < j \le r} \mathbb{E}_k[|A_{k+i} A_{k+j}|].$$

Put $\tilde{s} = s \wedge 2$. The Hölder inequality gives

$$\mathbb{E}_k[|A_{k+i} A_{k+j}|] \le \{\mathbb{E}_k[|A_{k+i}|^{\tilde{s}}]\}^{1/\tilde{s}} \{\mathbb{E}_k[|A_{k+j}|^{\tilde{s}/(\tilde{s}-1)}]\}^{(\tilde{s}-1)/\tilde{s}}.$$

Observing that $\frac{\tilde{s}}{\tilde{s}-1} - \tilde{s} \ge 0$ and using (C-1), we have, for all $\nu \in [0, \nu_1]$, $\mathbb{E}_k[|A_{k+j}|^{\tilde{s}/(\tilde{s}-1)}] \le \mathbb{E}_k[|A_{k+j}|^{\tilde{s}}]\nu^{-\tilde{s}/(\tilde{s}-1)+\tilde{s}}$, showing

$$\mathbb{E}_k[|S_k|] \le M_1\nu^{\tilde{s}} \left\{ \sup_{1 \le i \le r} \mathbb{E}_k[|A_{k+i}|^{\tilde{s}}] \right\}^{1/\tilde{s}} \left\{ \sup_{1 < j \le r} \mathbb{E}_k[|A_{k+j}|^{\tilde{s}}] \right\}^{(\tilde{s}-1)/\tilde{s}}$$

$$\le M_1\nu^{\tilde{s}} \sup_{1 \le i \le r} \mathbb{E}_k[|A_{k+i}|^{\tilde{s}}] \le M_1\nu^{\tilde{s}} \left( \sup_{1 \le i \le r} \mathbb{E}_k[|A_{k+i}|^s] \right)^{\tilde{s}/s}.$$

The proof of (100) then follows by bounding the above sup by a sum and by applying (C-3). □

Define $N_n := \sum_{l=0}^{[n/r]} \mathbf{I}(\phi_{lr} \leq R_1)$. We have the following.

Lemma A.2. *There exist $\alpha_0 > 0$ and $\gamma > 0$ such that, for all $n \geq 0$,*

$$\mathbb{E}_0[e^{-\alpha_0 N_n}] \leq e^{-\gamma n \alpha_0} V_0.$$

Proof.   Observe that, using (C-4), for all $k \in \mathbb{N}$,

$$\mathbb{E}_{kr}[V_{(k+1)r}] \leq V_{kr}[\lambda \mathbf{I}(\phi_{kr} > R_1) + B\mathbf{I}(\phi_{kr} \leq R_1)]$$

$$= V_{kr} \lambda^{\mathbf{I}(\phi_{kr} > R_1)} B^{\mathbf{I}(\phi_{kr} \leq R_1)}.$$

Let $\{W_{kr}, k \in \mathbb{N}\}$ be the process defined by

$$W_0 := V_0 \quad \text{and} \quad W_{kr} = \left(\frac{1}{\lambda}\right)^k \left(\frac{\lambda}{B}\right)^{N_{(k-1)r}} V_{kr}, \qquad k \geq 1.$$

Since $N_{kr}$ is $\mathcal{F}_{kr}$-measurable, we obtain, for all $k \in \mathbb{N}$,

$$\mathbb{E}_{kr}[W_{(k+1)r}] = \left(\frac{1}{\lambda}\right)^{k+1} \left(\frac{\lambda}{B}\right)^{N_{kr}} \mathbb{E}_{kr}[V_{(k+1)r}]$$

$$\leq \left(\frac{1}{\lambda}\right)^k \left(\frac{\lambda}{B}\right)^{N_{kr} - \mathbf{I}(\phi_{kr} \leq R_1)} V_{kr} = W_{kr}.$$

Hence, by induction $\mathbb{E}_0[W_{kr}] \leq W_0 = V_0$ and, since $V_k \geq 1$, we get

$$\mathbb{E}_0\left[\left(\frac{\lambda}{B}\right)^{N_{kr}}\right] \leq \lambda^{k+1} \mathbb{E}_0\left[\left(\frac{1}{\lambda}\right)^{k+1} \left(\frac{\lambda}{B}\right)^{N_{kr}} V_{(k+1)r}\right]$$

$$\leq \lambda^{k+1} \mathbb{E}_0[W_{(k+1)r}] \leq \lambda^{k+1} V_0.$$

Noting that $N_{kr+q} = N_{kr}$ for all $q = 0, 1 \dots, r-1$, the proof follows. □

We now turn back to the proof of Theorem 16. From (C-1), for all $k \in \mathbb{N}$ and $\nu \in [0, \nu_1]$,

$$\left|\prod_{l=kr+1}^{(k+1)r}(I - \nu A_l)\right|^p$$

(101)
$$\leq e^{-(\delta/2)\nu \mathbf{I}(\phi_{kr} \leq R_1)} \left\{\left|\prod_{l=kr+1}^{(k+1)r}(I - \nu A_l)\right|^p \right.$$

$$\left. \times e^{(\delta/2)\nu} \mathbf{I}(\phi_{kr} \leq R_1) + \mathbf{I}(\phi_{kr} > R_1)\right\},$$

where $\delta$ is defined in Lemma A.1. Let $n = mr + t$, where $m \in \mathbb{N}$ and $t = 0, 1, \ldots, r - 1$. Equation (101) and the Cauchy–Schwarz inequality show that

$$\mathbb{E}_0 \left| \prod_{l=1}^{n} (I - \nu A_l) \right|^p \leq \pi_1^{1/2} \pi_2^{1/2}, \qquad \text{where } \pi_1 = \mathbb{E}_0 [e^{-\delta \nu N_n}]$$

and

$$\pi_2 = \mathbb{E}_0 \left[ \prod_{k=0}^{m-1} \left\{ \left| \prod_{l=kr+1}^{(k+1)r} (I - \nu A_l) \right|^{2p} e^{\delta \nu} \mathbf{I}(\phi_{kr} \leq R_1) + \mathbf{I}(\phi_{kr} > R_1) \right\} \right].$$

Let $U_0 = 1$ and recursively define $U_{k+1}$ for $k = 0, 1, \ldots,$

$$U_{k+1} := \left\{ \left| \prod_{l=kr+1}^{(k+1)r} (I - \nu A_l) \right|^{2p} e^{\delta \nu} \mathbf{I}(\phi_{kr} \leq R_1) + \mathbf{I}(\phi_{kr} > R_1) \right\} U_k.$$

Applying Lemma A.1 with $a = 2p$, we obtain that $(U_k, \mathcal{F}_{kr})$ is a supermartingale. Consequently, $\pi_2 \leq 1$. Lemma A.2 and Jensen's inequality show that, for all $\nu \in [0, \alpha_0/\delta]$,

$$\mathbb{E}_0 [e^{-\delta \nu N_n}] \leq (\mathbb{E}_0 [e^{-\alpha_0 N_n}])^{\delta \nu / \alpha_0} \leq e^{-\gamma \delta \nu n} V_0,$$

which concludes the proof. $\square$

## APPENDIX B

**Burkholder inequalities for the TVAR process.** Throughout this section we let $\beta \in (0, 1]$, $L > 0$, $0 < \rho < 1$, $0 < \sigma_- \leq \sigma_+$ and we set $\mathcal{C} := \mathcal{C}(\beta, L, \rho, \sigma_-, \sigma_+)$. We further let $\tau \in (\rho, 1)$. The following lemma is adapted from [9], Proposition 4.

LEMMA B.1. *Let $(\Omega, \mathcal{F}, \mathbb{P}, \{\mathcal{F}_n; n \in \mathbb{N}\})$ be a filtered space. Let $p \geq 2$, $p_1, p_2 \in [1, \infty]$ such that $p_1^{-1} + p_2^{-1} = 2p^{-1}$ and let $\{Z_n; n \in \mathbb{N}\}$ be an adapted sequence such that $\mathbb{E}[Z_n] = 0$ and $\|Z_n\|_p < \infty$ for all $n \in \mathbb{N}$. Then*

$$(102) \qquad \left\| \sum_{i=1}^{n} Z_i \right\|_p \leq \left( 2p \sum_{i=1}^{n} \|Z_i\|_{p_1} \sum_{j=i}^{n} \|\mathbb{E} Z_j \mathcal{F}_i\|_{p_2} \right)^{1/2}.$$

PROPOSITION B.2. *Assume* (A1) *with $q \geq 2$ and let $p \geq 0$ be such that $2(p + 1) \leq q$. Then there exists $M > 0$ such that, for all $(\boldsymbol{\theta}, \sigma) \in \mathcal{C}$, $1 \leq s \leq t \leq n$ and sequences $\{\phi_i\}_{s \leq i \leq t}$ in $\mathrm{Li}(\mathbb{R}^d, 1, \mathbb{R}; p)$,*

$$\left\| \sum_{i=s}^{t} (\phi_i(\mathbf{X}_{i,n}) - \mathbb{E}_{\boldsymbol{\theta}, \sigma}[\phi_i(\mathbf{X}_{i,n})]) \right\|_{q/(p+1), \boldsymbol{\theta}, \sigma}^2 \leq M \sup_{i \in \{s, \ldots, t\}} |\phi_i|_{\mathrm{Li}(p)} \sum_{i=s}^{t} |\phi_i|_{\mathrm{Li}(p)}.$$

Proof. Let us apply Lemma B.1 with $Z_i = \phi_i(\mathbf{X}_{i,n}) - \mathbb{E}_{\boldsymbol{\theta},\sigma}[\phi_i(\mathbf{X}_{i,n})]$ and $p_1 = p_2 = q/(p+1)$. From the $L^q$ stability, we see that $\|Z_i\|_{q/(p+1),\boldsymbol{\theta},\sigma} \leq M|\phi_i|_{\mathrm{Li}(p)}$. It now remains to bound $\|\mathbb{E}_{\boldsymbol{\theta},\sigma}^{\mathcal{F}_{i,n}}[Z_k]\|_{q/(p+1),\boldsymbol{\theta},\sigma}$. Using the exponential stability and the $L^q$ stability, Proposition 11 shows that, for all $s \leq i \leq k \leq t$, $(\boldsymbol{\theta},\sigma) \in \mathcal{C}$, $\|\mathbb{E}_{\boldsymbol{\theta},\sigma}^{\mathcal{F}_{i,n}}[Z_k]\|_{q/(p+1),\boldsymbol{\theta},\sigma} \leq M\tau^{k-i}|\phi_k|_{\mathrm{Li}(p)}$. The proof follows. □

Another application of Lemma B.1 is the following result.

Proposition B.3. *Assume that $q \geq 5$ and let $p, r \geq 0$ be such that $u := 2(p+r) + 5 \leq q$. There exists a constant $M$ such that, for all $(\boldsymbol{\theta},\sigma) \in \mathcal{C}$, $1 \leq s \leq t \leq n$ and sequences $\{\gamma_{i,j}\}_{s \leq i < j \leq t}$ and $\{\phi_i\}_{s < i \leq t}$, respectively in $\mathrm{Li}(\mathbb{R}, 1, \mathbb{R}; p)$ and $\mathrm{Li}(\mathbb{R}, 1, \mathbb{R}; r)$,*

$$\left\| \sum_{s \leq i < j \leq t} \gamma_{i,j}(\mathbf{X}_{i,n})\sigma_{i+1,n}\varepsilon_{i+1,n}(\phi_j(\mathbf{X}_{j,n}) - \mathbb{E}_{\boldsymbol{\theta},\sigma}[\phi_j(\mathbf{X}_{j,n})]) \right\|_{2q/u,\boldsymbol{\theta},\sigma}$$

$$\begin{aligned}
(103) \quad &\leq M\Bigg\{ \sup_{s \leq i < j \leq t} |\gamma_{i,j}|_{\mathrm{Li}(p)} \sum_{i=s+1}^{t} |\phi_i|_{\mathrm{Li}(r)} \\
&\quad + \left( \sup_{s < i \leq t} |\phi_i|_{\mathrm{Li}(r)} \sum_{i=s+1}^{t} |\phi_i|_{\mathrm{Li}(r)} \right)^{1/2} \sup_{s < j \leq t} \left( \sum_{i=s}^{j} |\gamma_{i,j}|_{\mathrm{Li}(p)}^2 \right)^{1/2} \Bigg\}.
\end{aligned}$$

Proof. Let $\zeta_{i,j} := \gamma_{i-1,j}(\mathbf{X}_{i-1,n})\sigma_{i,n}\varepsilon_{i,n}$ and $U_j := \phi_j(\mathbf{X}_{j,n}) - \mathbb{E}_{\boldsymbol{\theta},\sigma}[\phi_j(\mathbf{X}_{j,n})]$ for all $s < i \leq j \leq t$. For all $s < i \leq j \leq t$, $U_i$ and $\zeta_{i,j}$ are $\mathcal{F}_{i,n}$-measurable. Throughout the proof, we denote by $M$ some constant independent of $s$, $t$, $n$ and $(\boldsymbol{\theta},\sigma) \in \mathcal{C}$. From uniform $L^q$ stability, for all $s < i \leq j \leq t$,

$$(104) \quad \|\zeta_{i,j}\|_{q/(p+1),\boldsymbol{\theta},\sigma} \leq M|\gamma_{i-1,j}|_{\mathrm{Li}(p)} \quad \text{and} \quad \|U_j\|_{q/(r+1),\boldsymbol{\theta},\sigma} \leq M|\phi_j|_{\mathrm{Li}(r)}.$$

Denote $Z_j := U_j \sum_{i=s+1}^{j} \zeta_{i,j}$. The LHS of (103) then reads

$$\begin{aligned}
(105) \quad &\left\| \sum_{s < i \leq j \leq t} \zeta_{i,j} U_j \right\|_{2q/u,\boldsymbol{\theta},\sigma} \\
&\leq \sum_{s < i \leq j \leq t} |\mathbb{E}_{\boldsymbol{\theta},\sigma}[\zeta_{i,j} U_j]| + \left\| \sum_{j=s+1}^{t} (Z_j - \mathbb{E}_{\boldsymbol{\theta},\sigma}[Z_j]) \right\|_{2q/u,\boldsymbol{\theta},\sigma}.
\end{aligned}$$

By Proposition 11, for all $s < i \leq j < t$,

$$(106) \quad \|\mathbb{E}_{\boldsymbol{\theta},\sigma}^{\mathcal{F}_{i,n}}[U_j]\|_{q/(r+1),\boldsymbol{\theta},\sigma} \leq M\tau^{j-i}|\phi_j|_{\mathrm{Li}(r)}.$$

Using (104) and (106), for all $s < i \le j \le t$,

$$|\mathbb{E}_{\boldsymbol{\theta},\sigma}[\zeta_{i,j} U_j]| = |\mathbb{E}_{\boldsymbol{\theta},\sigma}[\zeta_{i,j} \mathbb{E}_{\boldsymbol{\theta},\sigma}^{\mathcal{F}_{i,n}}[U_j]]|$$

$$\le \|\zeta_{i,j}\|_{2,\boldsymbol{\theta},\sigma} \|\mathbb{E}_{\boldsymbol{\theta},\sigma}^{\mathcal{F}_{i,n}}[U_j]\|_{2,\boldsymbol{\theta},\sigma}$$

$$\le M\tau^{j-i} |\gamma_{i-1,j}|_{\mathrm{Li}(p)} |\phi_j|_{\mathrm{Li}(r)}.$$

It then follows that

$$\sum_{s < i \le j \le t} |\mathbb{E}_{\boldsymbol{\theta},\sigma}[\zeta_{i,j} U_j]| \le M\gamma_{\infty,\infty}^{\dagger} \sum_{j=s+1}^{t} |\phi_j|_{\mathrm{Li}(r)} \sum_{i=s+1}^{j} \tau^{j-i} \le M\gamma_{\infty,\infty}^{\dagger} \phi_1^{\dagger},$$

where $\gamma_{\infty,\infty}^{\dagger} := \sup_{s \le i \le j \le t} |\gamma_{i,j}|_{\mathrm{Li}(p)}$ and $\phi_1^{\dagger} := \sum_{i=s+1}^{t} |\phi_i|_{\mathrm{Li}(r)}$; applying this bound to the RHS of (105), we obtain the first term of (103).

We now bound the second term in the RHS of (105). Applying Lemma B.1 with $p = q/u$, $p_1 = q/(p+r+2)$ and $p_2 = q/(p+r+3)$, this terms satisfies

$$\left\|\sum_{j=s+1}^{t} \overline{Z_j}\right\|_{2q/u,\boldsymbol{\theta},\sigma}$$

$$\le \left(\frac{4q}{u} \sum_{j=s+1}^{t} \|\overline{Z_j}\|_{q/(p+r+2),\boldsymbol{\theta},\sigma} \sum_{k=j}^{t} \|\mathbb{E}_{\boldsymbol{\theta},\sigma}^{\mathcal{F}_{j,n}}[\overline{Z_k}]\|_{q/(p+r+3),\boldsymbol{\theta},\sigma}\right)^{1/2},$$

where $\overline{Z_i} := Z_i - \mathbb{E}_{\boldsymbol{\theta},\sigma}[Z_i]$. Applying the Hölder inequality, (104) and the Burkholder inequality (see [14], Theorem 2.12) shows that, for all $s < j \le t$,

$$\|\overline{Z_j}\|_{q/(p+r+2),\boldsymbol{\theta},\sigma} \le M|\phi_j|_{\mathrm{Li}(r)} \gamma_{2,\infty}^{\dagger}$$

with

$$\gamma_{2,\infty}^{\dagger} := \sup_{s < j \le t} \left[\sum_{i=s+1}^{j} |\gamma_{i-1,j}|_{\mathrm{Li}(p)}^2\right]^{1/2}.$$

From the two last displays, we see that, in order to obtain the second term in the RHS of (103) and thus conclude the proof, it is now sufficient to show that, for all $s < j \le t$,

$$\sum_{k=j}^{t} \|\mathbb{E}_{\boldsymbol{\theta},\sigma}^{\mathcal{F}_j}[\overline{Z_k}]\|_{q/(p+r+3),\boldsymbol{\theta},\sigma} \le M\phi_{\infty}^{\dagger} \gamma_{2,\infty}^{\dagger},$$

where $\phi_{\infty}^{\dagger} := \sup_{s < i \le t} |\phi_i|_{\mathrm{Li}(r)}$. In fact, below we bound the LHS of the previous equation by $A_j + B_j$ and show this inequality successively for $A_j$ and $B_j$. Denoting $\overline{\zeta_{i,k} U_k} := \zeta_{i,k} U_k - \mathbb{E}_{\boldsymbol{\theta},\sigma}[\zeta_{i,k} U_k]$, we have

$$\sum_{k=j}^{t} \|\mathbb{E}_{\boldsymbol{\theta},\sigma}^{\mathcal{F}_j}[\overline{Z_k}]\|_{q/(p+r+3),\boldsymbol{\theta},\sigma} \le A_j + B_j,$$

where

$$A_j := \sum_{k=j}^{t} \left\| \sum_{i=s+1}^{j} \mathbb{E}_{\boldsymbol{\theta},\sigma}^{\mathcal{F}_{j,n}} [\overline{\zeta_{i,k} U_k}] \right\|_{q/(p+r+3),\boldsymbol{\theta},\sigma},$$

$$B_j := \sum_{k=j}^{t} \sum_{i=j+1}^{k} \| \mathbb{E}_{\boldsymbol{\theta},\sigma}^{\mathcal{F}_{j,n}} [\overline{\zeta_{i,k} U_k}] \|_{q/(p+r+3),\boldsymbol{\theta},\sigma}.$$

The bound on $A_j$ is obtained as follows. The centering term in the definition in $A_j$ may be forgotten by multiplying the leading term by a factor 2. Then we use the fact that $\zeta_{i,k} \in \mathcal{F}_{j,n}$ for all $i \leq j$ and all $k \geq i$ with the Hölder inequality, and finally apply (106) and (104) with the Burkholder inequality for martingale sequences to the obtained norms. These three steps read

$$A_j \leq 2 \sum_{k=j}^{t} \left\| \sum_{i=s+1}^{j} \mathbb{E}_{\boldsymbol{\theta},\sigma}^{\mathcal{F}_{j,n}} [\zeta_{i,k} U_k] \right\|_{q/(p+r+2),\boldsymbol{\theta},\sigma}$$

$$\leq 2 \sum_{k=j}^{t} \| \mathbb{E}_{\boldsymbol{\theta},\sigma}^{\mathcal{F}_{j,n}} [U_k] \|_{q/(r+1),\boldsymbol{\theta},\sigma} \left\| \sum_{i=s+1}^{j} \zeta_{i,k} \right\|_{q/(p+1),\boldsymbol{\theta},\sigma}$$

$$\leq M \sum_{k=j}^{t} \tau^{k-j} |\phi_k|_{\mathrm{Li}(r)} \left( \sum_{i=s+1}^{j} |\gamma_{i-1,k}|_{\mathrm{Li}(p)}^2 \right)^{1/2}$$

$$\leq M \phi_\infty^\dagger \gamma_{2,\infty}^\dagger.$$

It remains to show a similar inequality for $B_j$. From (104) and (106), for all $s < i \leq j \leq k \leq t$,

$$\| \mathbb{E}_{\boldsymbol{\theta},\sigma}^{\mathcal{F}_{i,n}} [\overline{\zeta_{j,k} U_k}] \|_{q/(p+r+2),\boldsymbol{\theta},\sigma}$$

$$\leq \| \mathbb{E}_{\boldsymbol{\theta},\sigma}^{\mathcal{F}_{j,n}} [\overline{\zeta_{j,k} U_k}] \|_{q/(p+r+2),\boldsymbol{\theta},\sigma}$$

(107)

$$\leq 2 \| \zeta_{j,k} \mathbb{E}_{\boldsymbol{\theta},\sigma}^{\mathcal{F}_{j,n}} [U_k] \|_{q/(p+r+2),\boldsymbol{\theta},\sigma}$$

$$\leq 2 \| \zeta_{j,k} \|_{q/(p+1),\boldsymbol{\theta},\sigma} \| \mathbb{E}_{\boldsymbol{\theta},\sigma}^{\mathcal{F}_{j,n}} [U_k] \|_{q/(r+1),\boldsymbol{\theta},\sigma}$$

$$\leq M \tau^{k-j} |\gamma_{j-1,k}|_{\mathrm{Li}(p)} |\phi_k|_{\mathrm{Li}(r)}.$$

This bound is in fact useful only when $k-j$ is large. We now derive another bound for the same quantity useful when $j-i$ is large. Since for all $i < j \leq k$, $\mathbb{E}_{\boldsymbol{\theta},\sigma}[\zeta_{j,k}] = \mathbb{E}_{\boldsymbol{\theta},\sigma}^{\mathcal{F}_{i,n}} [\zeta_{j,k}] = 0$, we have

$$\mathbb{E}_{\boldsymbol{\theta},\sigma}^{\mathcal{F}_{i,n}} [\overline{\zeta_{j,k} U_k}] = \mathbb{E}_{\boldsymbol{\theta},\sigma}^{\mathcal{F}_{i,n}} [\phi_k(\mathbf{X}_{k,n}) \gamma_{j-1,k}(\mathbf{X}_{j-1,n}) \sigma_{j,n} \varepsilon_{j,n}]$$

$$- \mathbb{E}_{\boldsymbol{\theta},\sigma} [\phi_k(\mathbf{X}_{k,n}) \gamma_{j-1,k}(\mathbf{X}_{j-1,n}) \sigma_{j,n} \varepsilon_{j,n}].$$

Note that $\sigma_{j,n}\varepsilon_{j,n} = X_{j,n} - \boldsymbol{\theta}_{j-1,n}^T \mathbf{X}_{j-1,n}$ is linear in $(\mathbf{X}_{j-1,n}, \mathbf{X}_{j,n})$. By Lemma C.2 and since $\sup_{j,n} |\boldsymbol{\theta}_{j-1,n}| < \infty$, $\gamma_{j-1,k}(\mathbf{X}_{j-1,n})\phi_k(\mathbf{X}_{k,n})\sigma_{j,n}\varepsilon_{j,n}$, as a mapping of $(\mathbf{X}_{j-1,n}, \mathbf{X}_{j,n}, \mathbf{X}_{k,n})$, belongs to $\mathrm{Li}(\mathbb{R}^d, 3, \mathbb{R}; p + r + 2)$ and its norm is bounded by $M|\gamma_{j-1,k}|_{\mathrm{Li}(p)}|\phi_k|_{\mathrm{Li}(r)}$. Hence, applying Proposition 11 and $L^q$ stability to the RHS of the previous display gives, for all $i < j \le k$,

$$\|\mathbb{E}_{\boldsymbol{\theta},\sigma}^{\mathcal{F}_{i,n}}[\overline{\zeta_{j,k}U_k}]\|_{q/(p+r+3),\boldsymbol{\theta},\sigma} \le M\tau^{j-i}|\gamma_{j-1,k}|_{\mathrm{Li}(p)}|\phi_k|_{\mathrm{Li}(r)}.$$

Combining with (107), we get

$$\|\mathbb{E}_{\boldsymbol{\theta},\sigma}^{\mathcal{F}_{i,n}}[\overline{\zeta_{j,k}U_k}]\|_{q/(p+r+3),\boldsymbol{\theta},\sigma} \le M\tau^{(j-i)\vee(k-j)}|\gamma_{j-1,k}|_{\mathrm{Li}(p)}|\phi_k|_{\mathrm{Li}(r)}.$$

Applying these bounds to the definition of $B_i$, since for all $i \le t$

$$\sum_{i \le j \le k \le t} \tau^{(j-i)\vee(k-j)} \le 2/(\tau(1-\tau)(1-\sqrt{\tau})),$$

we finally obtain $B_i \le M\gamma_{\infty,\infty}^\dagger \phi_\infty^\dagger \le M\gamma_{2,\infty}^\dagger \phi_\infty^\dagger$, which yields the proof. $\quad\square$

## APPENDIX C: TECHNICAL LEMMAS

LEMMA C.1. *Let $A$ be a positive semi-definite symmetric matrix and let $I$ denote the identity matrix with the same size as $A$. Then $|I - A| \le 1 - \lambda_{\min}(A) + 2|A|\mathbf{I}(|A| > 1)$.*

PROOF. Since $|\cdot|$ denotes the quadratic operator norm, we have $|I - A| = \max(1 - \lambda_{\min}(A), \lambda_{\max}(A) - 1)$. If $1 - \lambda_{\min}(A) \ge \lambda_{\max}(A) - 1$, the claimed inequality is trivially true. Since $A$ is positive semi-definite, $\lambda_{\max}(A) = |A|$. If $1 - \lambda_{\min}(A) < \lambda_{\max}(A) - 1$, $|I - A| = |A| - 1$. In addition, in this case, we necessarily have $|A| > 1$. Hence, the right-hand side of the claimed inequality in this case reads $1 - \lambda_{\min}(A) + 2|A| = 1 + |A| + \lambda_{\max}(A) - \lambda_{\min}(A) \ge |A| - 1$. $\square$

LEMMA C.2. *Let $(\mathsf{E}, |\cdot|_{\mathsf{E}})$ and $(\mathsf{F}, |\cdot|_{\mathsf{F}})$ be two normed spaces.*

(1) *Let $(\mathsf{G}, |\cdot|_{\mathsf{G}})$ be a normed space. For any $p_1, p_2 \ge 0$, there exists $C > 0$ such that, for all $\phi \in \mathrm{Li}(\mathsf{G}, 1, \mathsf{F}; p_1)$ and $\psi \in \mathrm{Li}(\mathsf{E}, m, \mathsf{G}; p_2)$,*

$$|\phi \circ \psi|_{\mathrm{Li}(p_1 p_2 + p_1 + p_2)} \le C|\phi|_{\mathrm{Li}(p_1)}(1 + |\psi|_{\mathrm{Li}(p_2)}^{p_1+1}).$$

(2) *Let $(\mathsf{G}, |\cdot|_{\mathsf{G}})$ be a normed algebra. For any $p_1, p_2 \ge 0$ and any integers $m_1, m_2 \ge 1$, there exists $C > 0$ such that, for all $\phi \in \mathrm{Li}(\mathsf{E}, m_1, \mathsf{G}; p_1)$ and $\psi \in \mathrm{Li}(\mathsf{F}, m_2, \mathsf{G}; p_2)$,*

$$|\phi\psi|_{\mathrm{Li}(p_1 + p_2 + 1)} \le C|\phi|_{\mathrm{Li}(p_1)}|\psi|_{\mathrm{Li}(p_2)}.$$

LEMMA C.3. *Let $\beta \geq 0$ and $\nu \in (0, 1)$. Then there exist constants $C_1, C_2$ depending only on $\beta$ such that*

$$\sup_{t>0} t^\beta (1-\nu)^t \leq C_1 \nu^{-\beta},$$

$$\sum_{s=1}^\infty (1-\nu)^s s^\beta \leq C_2 \nu^{-(1+\beta)},$$

*with the convention $0^0 = 1$. Assume now that $\beta > 0$. Then, as $\nu \downarrow 0$,*

(108)    $$\sum_{s=0}^\infty (I-\nu)^s ((s+1)^\beta - s^\beta) = \Gamma(\beta+1)\nu^{-\beta}(1 + O(\nu)),$$

*where $\Gamma$ is the Gamma function.*

PROOF. The result is trivial for $\beta = 0$, so we assume $\beta > 0$. A straightforward computation shows that $\sup_{t \geq 0} t^\beta (1-\nu)^t$ is attained at $t = t_0 := -\beta / \log(1-\nu)$. Since $\log(1-\nu)$ is bounded above by $-\nu$, the first bound is obtained. The second bound is obtained by bounding $(1-\nu)^s s^\beta$ by this sup for $s < t_0 + 1$ and by bounding the remainder of the sum (whose terms are decreasing) by

$$\int_0^\infty (1-\nu)^s s^\beta \, ds = \frac{1}{(-\log(1-\nu))^{\beta+1}} \int_0^\infty e^{-s} s^\beta \, ds,$$

which easily yields the first bound. For all $\beta > 0$ and $\nu \in (0, 1)$,

$$S := \sum_{s=0}^\infty (1-\nu)^s ((s+1)^\beta - s^\beta) = \beta \sum_{s=0}^\infty (1-\nu)^s \int_s^{s+1} t^{\beta-1} \, dt.$$

The proof of (108) follows by then writing

$$(1-\nu)S \leq \beta \int_0^\infty (1-\nu)^t t^{\beta-1} \, dt$$

$$= \Gamma(\beta+1)(-\log(1-\nu))^{-\beta} \leq S. \qquad \square$$

E. MOULINES
F. ROUEFF
GET/TÉLÉCOM PARIS
CNRS LTCI
46 RUE BARRAULT
75634 PARIS CEDEX 13
FRANCE
E-MAIL: moulines@tsi.enst.fr
        roueff@tsi.enst.fr

P. PRIOURET
LABORATOIRE DE PROBABILITÉS
UNIVERSITÉ PARIS VI
4 PLACE JUSSIEU
75252 PARIS CEDEX 05
FRANCE
E-MAIL: priouret@ccr.jussieu.fr